\documentclass[preprint,12pt]{elsarticle}

\usepackage{amssymb}
\usepackage{amsmath}
\usepackage{graphicx} 
\usepackage{mathtools}
\usepackage{tikz-cd}
\usepackage{amsthm}

\usepackage{lineno}

\newcommand{\p}{\mathbb{P}}

\newcommand{\A}{\mathbb{A}}
\newcommand{\TEmb}{\mathbb{T}}
\newcommand{\field}{\mathrm K}
\newcommand{\T}{\mathrm T}
\newcommand{\N}{\mathrm N}
\renewcommand{\subset}{\subseteq}
\renewcommand{\phi}{\varphi}
\renewcommand{\setminus}{\smallsetminus}
\renewcommand{\O}{\mathrm{O}}

\DeclareMathOperator{\Aut}{Aut}
\DeclareMathOperator{\orb}{Orb}
\DeclareMathOperator{\stab}{Stab}
\DeclareMathOperator{\rk}{rk}
\DeclareMathOperator{\im}{im}
\DeclareMathOperator{\codim}{codim}
\DeclareMathOperator{\End}{End}
\DeclareMathOperator{\Sym}{Sym}

\DeclareMathOperator{\GL}{GL}

\DeclareMathOperator{\SO}{SO}
\DeclareMathOperator{\ch}{ch}

\DeclareMathOperator{\PO}{PO}

\newtheorem{theorem}{Theorem}[section]
\newtheorem{corollary}[theorem]{Corollary}
\newtheorem{lemma}[theorem]{Lemma}
\newtheorem{proposition}[theorem]{Proposition}
\newtheorem{definition}[theorem]{Definition}
\newtheorem{claim}{Claim}
\newtheorem{remark}[theorem]{Remark}

\journal{Advances in Applied Mathematics}

\begin{document}

\begin{frontmatter}


\author{Franquiz Caraballo Alba \corref{cor1}}
\ead{fc19b@fsu.edu}
\ead[url]{https://www.math.fsu.edu/~fcarabal/}
\cortext[cor1]{}
\affiliation{organization={Florida State University},
             addressline={1017 Academic Way},
             city={Tallahassee},
             postcode={32306},
             state={Florida},
             country={United States of America}
}

\title{Linear Orbits of Smooth Quadric Surfaces}

\begin{abstract}
The \textit{linear orbit} of a degree $d$ hypersurface in $\p^n$ is its orbit under the natural action of $\p \GL(n+1)$, in the projective space of dimension $N =\binom{n+d}{d} - 1$ parameterizing such hypersurfaces. This action restricted to a specific hypersurface $X$ extends to a rational map from the projectivization of the space of matrices to $\p^N$. The class of the graph of this map is the \textit{predegree polynomial} of its corresponding hypersurface. The objective of this paper is threefold. First, we formally define the predegree polynomial of a hypersurface in $\p^n$, introduced in the case of plane curves by Aluffi and Faber, and prove some results in the general case. A key result in the general setting is that a partial resolution of said rational map can contain enough information to compute the predegree polynomial of a hypersurface. Second, we compute the leading term of the predegree polynomial of a smooth quadric in $\p^n$ over an algebraically closed field with characteristic 0, and compute the other coefficients in the specific case $n = 3$. In analogy to Aluffi and Faber's work, the tool for computing this invariant is producing a (partial) resolution of the previously mentioned rational map which contains enough information to obtain the invariant. Third, we provide a complete resolution of the rational map in the case $n=3$, which in principle could be used to compute more refined invariants.
\end{abstract}

\begin{keyword}
    projective linear group, smooth quadric surfaces, algebraic geometry, intersection theory, resolution of singularities.

    \MSC[2020] 14N10 \sep 14L30 \sep 14C17.
\end{keyword}

\end{frontmatter}


\section{Introduction} \label{sec:introduction}
    Throughout this paper, $\field$ is an algebraically closed field of characteristic $0$. The group $\p \GL_\field(n+1)$ of projective linear transformations of $\p_\field^n$ acts on the space of hypersurfaces of a fixed degree. The study of the orbits of this action is a natural task; the orbit of a hypersurface depends subtly on its geometry, and invariants such as the degree of the closure of an orbit have a natural enumerative interpretation and are of interest in representation theory. The case $n = 1$ of orbits of point configurations in $\p^1$ was studied in \cite{aluffiFaberPoints}, while the considerably more involved case $n = 2$ was carried out by Aluffi and Faber in a sequence of papers (\cite{aluffiFaberCurves}, \cite{aluffiFaberAPP3}, \cite{aluffiFaberAPP}, \cite{aluffiFaberAPP2}) leading to the computation of the degree of the orbit closure of an arbitrary plane curve. For $n > 2$, very little is known. Tzigantchev (\cite{tzigantchev}) analyzed the case of plane configurations in $\p^3$, and using the same approach as Aluffi and Faber, recently Cazzador and Skauli (\cite{cazzadorSkauli}) gave a partial treatment of an intersection-theoretic computation of the degree of the orbit closure of a smooth cubic surface in $\p^3$. In \cite{numericalAPP}, Brunstega I Moncusi, Timme and Weinstein compute this degree numerically to be 96120. Deopurkar, Patel and Tseng recover this result in \cite{equivariantAPP} using techniques from equivariant geometry.

    In this paper we study the case of smooth quadric hypersurfaces in $\p^n$. The space of quadrics is a $N = \frac{(n+1)(n+2)}{2}-1$ dimensional projective space, and all smooth quadrics are projectively equivalent, so that the orbit of an arbitrary smooth quadric under the $\p \GL(n+1)$ action is dense in $\p^N$. Thus, the degree of the orbit closure of a smooth quadric is trivially 1. However, there are more refined invariants associated with these actions. Aluffi and Faber introduced in \cite{aluffiFaberAPP} the notion of ‘predegree polynomial’ of a hypersurface $X \subset \p^n$: the coefficient of the $i$th degree term of this polynomial is the number of $\p \GL(n+1)$ translates $\phi (X)$ containing $i$ general points, subject to imposing $n^2+2n-i$ general linear conditions on $\phi$. The precise definition of the predegree polynomial is in terms of the Chern character of an associated line bundle, or equivalently the multidegrees of a certain class in a product of projective spaces; cf.~\S\ref{sec:predegree_polynomial}. For instance, in the case of quadric hypersurfaces in $\p^n$, the leading coefficient of the polynomial equals the degree of the closure of the stabilizer of a smooth quadric in $\p^n$, that is, the degree of (the closure of) $\PO (n + 1)$ (cf.~Proposition \ref{prop:degreePO}). In \S\ref{sec:pp_quadric}, we obtain the following result:

    \begin{theorem}\label{thm:main}
        The predegree polynomial of a smooth quadric in $\p^3_\field$ is $1 + 2t + 4t^2 + 8t^3 + 16t^4 + 32t^5 + 64t^6 + 112t^7 + 140t^8 + 40t^9$.
    \end{theorem}

    This computes the degree of $\overline{\PO(4)} \subset \p^{15}$ to be $40$ and provides new enumerative results. For example, it shows that there are 140 translates $\phi(Q)$ of a smooth quadric containing 8 general points, with $\phi$ ranging in a general 8-dimensional linear space of projective transformations (for details see Corollary \ref{cor:interpretation} and Table \ref{table:interpretation}). We prove Theorem~\ref{thm:main} by producing a partial resolution of the indeterminacies of the rational map $\p^{15} \dashrightarrow \p^9$ extending the $\p \GL(4)$ action on a given quadric $Q$, and performing intersection-theoretic computations. This approach is analogous to the one employed by Aluffi and Faber in the case of configurations of points in $\p^1$ and of plane curves. The partial resolution is obtained by a sequence of blow-ups at smooth centers of the space $\p^{15}$ compactifying $\p \GL(4)$. In fact, while the partial resolution suffices for our enumerative purposes, we describe in §4 a sequence of blow-ups of $\p^{15}$ at smooth centers that completely resolve this rational map. This sequence could be used to compute other invariants associated with the $\p \GL(4)$ action and is of independent interest. We note that in the case of smooth cubic surfaces, a complete resolution of the corresponding rational map is not known. It would be interesting to use the partial resolution obtained in~\cite{cazzadorSkauli} to compute at least some of the coefficients of the predegree polynomial for a smooth cubic surface.

    The paper is organized as follows. In \S\ref{sec:background} we recall basic background, formulate precisely the question we answer for quadric surfaces in Theorem~\ref{thm:main} and prove results in the general setting. These results are useful in computing predegree polynomials for arbitrary degree and ambient dimension.
    
    In \S\ref{sec:pp_quadric_n} we address the $n$ dimensional case for quadrics, including a computation of the leading coefficient of the predegree polynomial of a smooth quadric in $\p^n$. The partial resolution needed for the proof of Theorem~\ref{thm:main} is obtained in \S \ref{sec:pp_quadric}, and we give details of the intersection-theoretic computation yielding the result. In \S\ref{sec:resolution} we describe the sequence of blow-ups needed for a complete resolution of the rational map. This essentially is a resolution of the singularities of the graph of the key rational map and could be used to compute other invariants of the graph, such as its multiplicity along different strata.

\section{Background} \label{sec:background}

    The group $\p \GL(n+1)$ acts on the space of hypersurfaces of $\p^n$ of fixed degree $d$. We are interested in the orbits of this action, from an enumerative point of view. The relevant question is `how many translates of a given hypersurface $X$ by matrices satisfying $n^2+2n-i$ general linear conditions contain $i$ general points?'. Here, each translate is counted with a multiplicity accounting for the number of automorphisms of $X$ induced by the matrices satisfying the given conditions. This multiplicity is $1$ if $i < \dim \orb(X)$ and if $i = \dim \orb(X)$, then it is $\deg(\overline{\stab(X)})$. For details, see Proposition \ref{prop:multiplicities}.
    
    The \textit{predegree polynomial} of a hypersurface is $\sum_i a_i t^i$, where $a_i$ is the number of these translates for a fixed value of $i$. For a rigorous definition and interpretation, see Definition~\ref{def:predegree_polynomial} and Corollary \ref{cor:interpretation}. In \cite{aluffiFaberAPP3}, the predegree polynomial is defined and computed in the case $n=2$, for arbitrary plane curves. The case $n=3$ where the surface $X$ is the union of planes is studied in \cite{tzigantchev}. In this paper, the focus is smooth quadrics in $\p^n$. In this section the relevant question is rigorously defined and some results about the general setting are proved.
    
    Let $A$ be an $n+1$ dimensional vector space over $\field$. Since $\Aut(\p A)$ is $\p \GL(A)$ and $\GL(A)$ embeds into $\End (A)$ as the complement of the zero locus of the determinant, $\Aut(\p A)$ embeds into the set of ``projective matrices'', $\p \End(A)$.
    
    For the rest of the paper, we abuse notation in the following sense. Given $\phi \in \p \End (A)$, let $\hat \phi \in \End(A)$ be a representative of $\phi$. We call the \textit{kernel} of $\phi$ the projectivization $\p \ker \hat{\phi}$ of the kernel of $\hat \phi$ inside $\p A$ and denote it by $\p \ker \phi$. The projective matrix $\phi$ determines a rational map $\p A \dashrightarrow \p A$ that is regular away from its kernel in the following sense: given $\hat p \in A$, a representative of $p \in \p A$ not in the kernel of $\phi$, we define the image of $p$ to be the element of $\p A$ represented by $\hat \phi (\hat v)$. We denote this rational map by $\phi$. Let $X$ be a subvariety of $\p A$ and let $\hat X$ be the affine cone over $X$ whose vertex is the origin $0$ in $A$. We call the closure of the image under the projectivization $A\setminus\{0\} \to \p A$ of $\hat \phi (\hat X \setminus \{0\})$ the \textit{image of $X$ under $\phi$} and denote it by $\overline{\phi(X)}$. The \textit{image of $\phi$} is $\im \phi := \overline{\phi(\p A)}$. Note that the projective matrix $\phi$ is invertible if and only if its kernel is empty or, equivalently, its image is $\p A$. Note that $\p \End(A)$ is no longer a ring or even a group, but instead a monoid, where the operation is inherited from the multiplication of $\End(A)$. The largest group in the monoid $\p \End(A)$ is precisely $\p \GL(A)$, and for $\phi, \psi \in \p \End(A)$, the map $\phi\psi$ is the composition $\phi \circ \psi$.

    \subsection{Point Conditions}\label{sec:point_conditions}
        Let $X$ be a degree $d$ hypersurface in $\p A$ and let $f \in \Sym^d(A^\vee)$ be a generator of its homogeneous ideal. Note that this choice is unique up to scalar, so choosing this hypersurface is equivalent to choosing a point in $\p \Sym^d(A^\vee)$. Furthermore, given $\phi \in \p \GL (A)$, the ideal generated by $f(\hat \phi(\_))$ in $\Sym^d(A^\vee)$ is independent of the choice of representative $\hat \phi$ of $\phi$. Let $X \circ \phi$ denote the hypersurface defined by $f(\hat \phi(\_))$.
        
        A straightforward check shows that for arbitrary elements $\phi,\psi$ of $\p \GL(A)$ $(X \circ \phi) \circ \psi = X \circ (\phi \circ \psi)$, making this a right action of $\p \GL(A)$ on $\p \Sym^d(A^\vee)$. For the rest of this section, we fix the hypersurface $X$, and thus, the class of the polynomial $f$ in $\p \Sym^d(A^\vee)$. The action of $\p \GL(A)$ defines a regular map $\p\GL(A) \to \p \Sym^d(A^\vee)$ by $\phi \mapsto X \cdot \phi$, which extends to a rational map on each compactification of $\p \GL(A)$. 

        \begin{definition} \label{def:alpha}
            The map $\alpha: \p \End(A) \dashrightarrow \p \Sym^d(A^\vee)$ is the rational map which extends the regular map $\p \GL(A) \to \p \Sym^d(A^\vee)$ by $\phi \mapsto X \cdot \phi$.
        \end{definition}
        
        Concretely, given $\phi \in \p \End(A)$, if $f(\phi(\_))$ is not identically zero, then $\alpha(\phi)$ is the hypersurface defined by $f(\phi(\_)) = 0$ and undefined otherwise. This map and its (partial) resolutions are central to our approach. In \ref{sec:resolution}, we resolve this map in the case $n=3$ through blowing up nonsingular subschemes and study the centers of the blowups in more detail. This requires a careful study of the ``indeterminacy scheme'' of $\alpha$; the subscheme of $\p \End(A)$ where $\alpha$ is undefined.

        Given a point $q\in \p A$, we obtain a hyperplane $H_q$ in $\p \Sym^d(A^\vee)$, consisting of the hypersurfaces $X$ containing $q$. This is a hyperplane since the condition $f(\hat q) = 0$ is linear in the entries of $\p \Sym^d(A^\vee)$.
        
        \begin{definition}[Point condition corresponding to a point in $\p A$]\label{def:pointCondition}
            The \textit{point condition corresponding to $q$}, which we denote by $P_q$, is the closure of the subspace of $\p \End(A)$ consisting of projective matrices whose image under $\alpha$ contains the given point $q \in \p A$.
        \end{definition}

       Since the set of hyperplanes $H_q$ spans the dual of $\p \Sym^d(A^\vee)$, the set of point conditions $P_q$ spans the linear system defining the rational map $\alpha$. Thus, the indeterminacy scheme of $\alpha$ is the (scheme-theoretic) intersection of the point conditions $P_q$. Therefore, $\phi \in \p \End(A)$ is in the support of the indeterminacy scheme of $\alpha$ if and only if $f(\phi(q)) = 0$ for all $q \in \p A$. In other words, $\phi$ is in the support of the indeterminacy scheme of $\alpha$ if and only if $\im \phi \subset X$. We record this and the fact that the closure of the image of $\alpha$ is the closure of the orbit of $X$ in the following proposition, stated without proof.
        
        \begin{proposition}\label{prop:base_scheme}
            The support of the indeterminacy scheme of $\alpha$ is $B_0 := \{\phi \in \p \End(A) | \im \phi \subset X\}$. Furthermore, the closure of the image 
            of $\alpha$ is the closure of the orbit of $X$ in $\p \Sym^d(A^\vee)$ under the action of $\p \GL(A)$.
        \end{proposition}

        For the rest of this subsection, it will be convenient to fix a choice of homogeneous coordinates $x_0,\dots,x_{n}$ on $\p A$.

        \begin{proposition}\label{prop:tangent_point_condition}
            Given $q \in \p A$ and $\phi \in P_q$ such that $q \notin \p \ker \phi$ and $X$ is nonsingular at $\phi (q)$, the embedded tangent space to $P_q$ at $\phi$ consists of the projective matrices $\psi \in \p \End(A)$ such that $q \in \p \ker \psi$ or $\psi(q)$ is in the embedded tangent space to $X$ at $\phi(q)$.
    	\end{proposition}
        \begin{proof}
             Let $\psi \in\p \End(A)$, $q \in \p A$ and $\phi \in P_q$ be such that $q \notin \p \ker \phi$. Let $\hat \psi $ be a representative of $\psi$, $\hat \phi$ be a representative of $\phi$ and $v$ a representative of $q$.
            
            With this notation, $\psi$ is in the embedded tangent space to $P_q$ at $\phi$ if and only if $\nabla s_q(\hat \phi) \cdot \hat \psi  = 0$, where $s_q = f(\_(v))$, $\nabla$ represents the gradient and $\_ \cdot \_$ is the standard dot product for our choice of basis of $A$. Then, using the chain rule, $\nabla s_q(\hat \phi) \cdot \hat \psi = \nabla f(\hat \phi v) \cdot \hat \psi v$. Thus, expanding the expression, $\nabla s_q(\hat \phi) \cdot \hat \psi  = 0$ if and only if:
            \begin{equation}\label{eqn:tangent_P_q}
                \sum_{i=0}^{n} \frac{\partial f}{\partial x_i} \bigg\vert_{\hat \phi v} (\hat \psi v)_i = 0
            \end{equation}
            where the subscripts indicate the relevant projection to $\field$. 
    
            Note that \eqref{eqn:tangent_P_q} holds if and only if $\hat \psi v$ is in the hyperplane dual to $\nabla f(\hat \phi v)$ or $\nabla f(\hat \phi v) = 0$.
    
            Since $\hat \phi$ is a representative of $\phi$ and $v$ is a representative of $q$, $\hat \phi v$ is a representative of $\phi(q)$. By assumption, $q \notin \p \ker \phi$, so $\hat \phi v \neq 0$. Furthermore, for $q' \in \p A$, $q'$ is in the embedded tangent space of $X$ at $\phi(q)$ if and only if for all representatives $w$ of $q'$ it is true that $\nabla f(\hat \phi v) \cdot w = 0$; that is, if and only if
            \begin{equation}\label{eqn:tangent_X}
               \sum_{i=0}^{n} \frac{\partial f}{\partial x_i} \bigg\vert_{\hat \phi v} (w)_i = 0. 
            \end{equation}
            Since by assumption $X$ is nonsingular at $\phi (q)$, $\nabla f(\hat \phi v) \neq 0$. Thus, $\psi$ is in the embedded tangent space to $\phi$ if and only if $\hat \psi v$ is in the hyperplane dual to $\nabla f(\hat \phi v)$.
    
            Assume that $\psi$ is in the embedded tangent space to $P_q$ at $\phi$. This means that $\nabla f(\hat \phi v) \cdot \hat \psi v = 0$, so \eqref{eqn:tangent_X} holds for $w = \hat \psi v$; that is, $\psi (q)$ is in the embedded tangent space to $X$ at $\phi (q)$.
    
            Now, assume that $\psi (q)$ is in the embedded tangent space to $X$ at $\phi (q)$. Then, $\nabla s_q(\hat \phi) \cdot \hat \psi  = \nabla f(\hat \phi v) \cdot \hat \psi v = 0$; that is, \eqref{eqn:tangent_P_q} holds, thus $\psi$ is in the embedded tangent space to $P_q$ at $\phi$.
        \end{proof}
    
        \begin{corollary}\label{cor:point_condition_singular}
            Let $q \in \p A$ and $\phi \in P_q$. If $P_q$ is singular at $\phi$, then $q \in \p \ker \phi$ or $X$ is singular at $\phi(q)$. 
        \end{corollary}
        \begin{proof}
            Let $q \in \p A$ and $\phi \in P_q$. Assume that $q \notin \p \ker \phi$ and that $P_q$ is singular at $\phi$. Let $v$ be a representative of $q$ and $\hat \phi$ be a representative of $\phi$. Then, by the discussion in Proposition~\ref{prop:tangent_point_condition}, this means that $\nabla s_q (\hat \phi) = 0$, thus the embedded tangent space to $P_q$ at $\phi$ is all of $\p \End(A)$.
    
            Let $q' \in \p A$ be arbitrary and $w$ be a representative of $q'$. Then, there exists some invertible matrix $\hat \psi  \in \GL(A)$ such that $\hat \psi v = w$. Since $P_q$ is singular at $\phi$, the class of $\hat \psi $ in $\p \End(A)$ is in the embedded tangent space to $P_q$ at $\phi$; that is, $\nabla s_q(\hat \phi) \cdot \hat \psi  = 0$. Then:
            \begin{align*}
                \nabla f(\hat \phi v) \cdot w & = \nabla f(\hat \phi v) \cdot \hat \psi v \\
                & = \nabla s_q(\hat \phi) \cdot \hat \psi  \\
                & = 0.
            \end{align*}
            Thus, $q'$ is in the embedded tangent space to $X$ at $\phi(q)$, so $X$ is singular at $\phi(q)$.
        \end{proof}

        \begin{remark} \label{rmk:point_condition_smooth}
            If $q \notin \p \ker \phi$ and $X$ is smooth at $\phi(q)$, then $P_q$ is smooth at $\phi$.
        \end{remark}
        
    \subsection{The Predegree Polynomial}\label{sec:predegree_polynomial}
    
        Consider the closure $\overline{\Gamma}$ of the graph of $\alpha$ in $\p \End(A) \times \p \Sym^d(A^\vee)$ and let $p_1$, resp., $p_2$ be the natural projections to $\p \End(A)$, resp., $\p \Sym^d(A^\vee)$.
        \begin{center}
            \begin{tikzcd}
                \overline{\Gamma} \arrow[d, "p_1", two heads] \arrow[rrd, "p_2"] &  &  \\
                \p \End(A) \arrow[rr, "\alpha", dashed] &  & \p \Sym^d(A^\vee)
            \end{tikzcd}
        \end{center}
        Note that $\im p_2$ is the closure of the orbit of $X$ under the action of $\p \GL(A)$, which we denote by $\overline{\orb(X)}$. Also note that $\dim \overline{\Gamma} = \dim \p \End(A) = n^2+2n$.
        We define the multidegrees of $\alpha$ as follows:
        \begin{equation}\label{eqn:pp_def1}
            a_i = \int_{\p \End(A) \times \p \Sym^d(A^\vee)}  p_1^* h_1^{n^2+2n-i} p_2^* h_2^{i} \cap [ \overline{\Gamma} ]
        \end{equation}
        where $h_1 = c_1(\mathcal{O}_{\p \End(A)}(1))$, $h_2 = c_1(\mathcal{O}_{\p \Sym^d(A^\vee)}(1))$ and $i$ ranges from $0$ to $n^2+2n$. Pushing forward by $p_2$ gives
        \begin{equation}\label{eqn:pp_def1_p2}
            a_i = \int_{\overline{\orb(X)}} h_2^i \cap p_{2*} \left(p_1^* h_1^{n^2+2n-i} \cap [\overline{\Gamma}] \right).
        \end{equation}
        Therefore, if $i > \dim \overline{\orb(X)}$, then $a_i = 0$. Pushing forward by $p_1$ gives
        \begin{equation}\label{eqn:pp_def1_p1}
            a_i = \int_{\p \End(A)} h_1^{n^2+2n-i} p_{1*}\left(p_2^* h_2^{i} \cap [ \overline{\Gamma} ]\right).
        \end{equation}

        Let $\pi:\widetilde{S} \to \p \End(A)$ be a proper birational map resolving the indeterminacies of $\alpha$; that is, equipped with a regular map $\widetilde{\alpha}: \widetilde{S} \to \p \Sym^d(A^\vee)$ such that the following diagram commutes:
        \begin{center}
            \begin{tikzcd}
                \widetilde{S} \arrow[d, "\pi", two heads] \arrow[rrd, "\widetilde{\alpha}"] &  &  \\
                \p \End(A) \arrow[rr, "\alpha", dashed] &  & \p \Sym^d(A^\vee)
            \end{tikzcd}
        \end{center}    
        Define $b_i$ by
        \begin{equation}\label{eqn:pp_def2}
            b_i = \int_{\p \End(A)} h_1^{n^2+2n-i} \cap \pi_* \left( \widetilde{\alpha}^* h_2^i \cap [\widetilde{S}] \right).
        \end{equation}
    
        Finally, consider the same resolution $\pi:\widetilde{S} \to \p \End(A)$ of $\alpha$ and let $m$ be the dimension of $\p \Sym^d(A^\vee)$. Let $q_1, \dots, q_{m} \in \p A$ be points in general position. Let $Y_i$ denote the intersection of $\widetilde{P_{q_1}}, \dots, \widetilde{P_{q_i}}$, where $\widetilde{P_q} = \widetilde{\alpha}^{-1}(H_q)$ and $j$ is the inclusion of $Y_i$ in $\widetilde{S}$. Define $c_i$ by         
        \begin{equation}\label{eqn:pp_def3}
            c_i = \int_{Y_i} j^* \pi^*h_1^{n^2+2n-i} \cap [Y_i].
        \end{equation}

        \begin{proposition}\label{prop:equivalence_def}
            For all $i$, $a_i= b_i = c_i$. In particular, $b_i$ and $c_i$ are independent of the choice of the resolution $\pi$.
        \end{proposition}
        \begin{proof}
            There exists a map $\eta:\widetilde{S} \to \overline{\Gamma}$ such that the following diagram commutes:
            \begin{center}
                \begin{tikzcd}
                    \widetilde{S} \arrow[rdd, "\pi"', two heads] \arrow[rrd, "\widetilde{\alpha}"] \arrow[rd, "\eta"] & & \\
                    & \overline{\Gamma} \arrow[r, "p_2", two heads] \arrow[d, "p_1"', two heads] & \p \Sym^d(A^\vee) \\
                    & \p \End(A) \arrow[ru, "\alpha", dashed] &
                \end{tikzcd}
            \end{center}
            Then, by functoriality, $\widetilde{\alpha}^* = \eta^* p_2^*$ and $\pi_* = p_{1*} \eta_*$. Thus,
            \begin{align*}
                b_i = & \int_{\p \End(A)} h_1^{n^2+2n-i} \cap \pi_* \left( \widetilde{\alpha}^* h_2^i \cap [\widetilde{S}] \right) \\
                = & \int_{\p \End(A)} h_1^{n^2+2n-i} \cap p_{1*} \eta_* \left( \eta^* p_2^* h_2^i \cap [\widetilde{S}] \right) \\
                = & \int_{\p \End(A)} h_1^{n^2+2n-i} \cap p_{1*} \left( p_2^* h_2^i \cap \eta_* [\widetilde{S}] \right) \\
                = & \int_{\p \End(A)} h_1^{n^2+2n-i} \cap p_{1*} \left( p_2^* h_2^i \cap [\overline{\Gamma}] \right) \\
                = & a_i,
            \end{align*}
            where the last line uses equation \eqref{eqn:pp_def1_p1}. Since $\widetilde{P_q} = \widetilde{\alpha}^{-1}(H_q)$ and $H_q$ is a hyperplane, $j_*[Y_i] = \widetilde{\alpha}^*h_2^i \cap [\widetilde{S}]$. Pushing forward the intersection-theoretic computation defining $c_i$ by $j$:
            
            \begin{align*}
                c_i = & \int_{Y_i} j^* \pi^*h_1^{n^2+2n-i} \cap [Y_i] \\
                = & \int_{\widetilde{S}} j_* \left(j^* \pi^*h_1^{n^2+2n-i} \cap [Y_i]\right) \\
                = & \int_{\widetilde{S}} \pi^*h_1^{n^2+2n-i} \cap j_* \left( [Y_i]\right) \\
                = & \int_{\widetilde{S}} \pi^*h_1^{n^2+2n-i} \cap \left( \widetilde{\alpha}^*h_2^i \cap [\widetilde{S}] \right) \\
                = & \int_{\p \End(A)} \pi_* \left(\pi^* h_1^{n^2+2n-i} \cap \left( \widetilde{\alpha}^*h_2^i \cap [\widetilde{S}] \right) \right) \\
                = & \int_{\p \End(A)} h_1^{n^2+2n-i} \cap \pi_* \left( \widetilde{\alpha}^*h_2^i \cap [\widetilde{S}] \right) \\
                = & b_i,
            \end{align*}
            as desired.
        \end{proof}

        \begin{definition}\label{def:predegree_polynomial}
            The `predegree polynomial' of a given hypersurface $X$ in $\p A$ is 
            \begin{equation*}
                \mathcal{P}_X(t) = \sum_{i=0}^{n^2+2n} a_i t^i
            \end{equation*}
            where $a_i$ is defined as in \eqref{eqn:pp_def1} above.
        \end{definition}

       In \cite{aluffiFaberPoints}, Aluffi and Faber defined the \textit{predegree} of the orbit closure of a hypersurface $X$ to be the degree of the orbit closure of $X$ times the degree of the stabilizer of $X$ because it was the natural output of the intersection-theoretic computation preliminary to the computation of the degree of the orbit closure. The predegree polynomial is called so because it is a direct generalization of the intersection-theoretic information. For an interpretation of the coefficients of the predegree polynomial, see Corollary \ref{cor:interpretation}.
        
        The definition of $c_i$ reflects the informal question posed at the beginning of this section. This comes from the fact that, given $q \in \p A$, $H_q$ is the hyperplane of $\p \Sym^d(A^\vee)$ consisting of the hypersurfaces $X$ such that $q \in X$. Then, $\widetilde{P_q} = \widetilde{\alpha}^{-1}(H_q)$ consists of the matrices sending $X$ to hypersurfaces containing $q$, so $Y_i$ consists of the matrices which send $X$ to a hypersurface containing the points $q_1, \dots, q_i$. Then, the degree of $Y_i$ is the number of matrices $\phi$ subject to $n^2+2n-i$ general linear conditions such that $\phi(X)$ contains $i$ general points.

         \begin{proposition}\label{prop:multiplicities}
            Let $L$ be a general dimension $i$ linear subspace of $\p \End(A)$ and let $X'$ be a translate of $X$ by an element of $\p \GL(A) \cap L$. If $i < \dim \orb(X)$, there exists a unique $\phi \in L \cap \p \GL(A)$ such that $X' = X \circ \phi$. If $i = \dim \orb(X)$, then there are $\deg \overline{\stab(X)}$ elements $\phi \in L$ such that $X'= X \circ \phi$.
        \end{proposition}
        \begin{proof}
            Assume that $i < \dim \orb(X)$ and, WLOG, assume that $X' = X$. Since $X$ is a translate of $X$ by an element of $L \cap \p \GL(A)$, there exists some $\phi \in L$ such that $X \circ \phi = X$. Therefore, $\phi \in \overline{\stab(X)} \cap L$. By the orbit stabilizer theorem, $\overline{\stab(X)}$ has codimension $\dim \orb(X)$ and thus cannot intersect a general $L$ at more than one point. Thus $\phi$ is unique.

            The second part of the statement follows from the definition of degree.
        \end{proof}

        The information captured by the predegree polynomial can be interpreted as follows:

        \begin{corollary}\label{cor:interpretation}
            Given a hypersurface $X \subset \p A$ with predegree polynomial $\mathcal{P}_X(t) = a_0 + a_1t + \dots + a_{n^2+2n}t^{n^2+2n}$:
            \begin{itemize}
                \item If $i < \dim \orb(X)$, then $a_i$ is the number of $L$-translates of $X$ containing $i$ general points, where $L$ is the intersection of $\p \GL(A)$ and a general $i$ dimensional linear subspace of $\p \End(A)$.
                \item If $i = \dim \orb(X)$, then $a_i$ is the number of translates of $X$ containing $\dim \orb(X)$ general points times the degree of the closure of the stabilizer of $X$. That is, it is the degree of the closure of the orbit of $X$ times the degree of the closure of the stabilizer of $X$.
                \item If $i > \dim \orb(X)$, then $a_i = 0$.
            \end{itemize}
        \end{corollary}

        The information of the predegree polynomial is not only equivalent to the information of the multidegrees of the map $\alpha$, but it is also equivalent to the information of the pullback of the Chern character of $\mathcal{O}_{\p \Sym^d(A^\vee)}(1)$. In fact
        
        \begin{remark}\label{rem:chern_character}
            For any proper birational map $\pi: \widetilde{S} \to \p \End(A)$ resolving the indeterminacies of $\alpha$,
            \begin{equation*}
                \sum_{i=0}^{n^2+2n} \frac{a_i}{i!}h_1^i \cap [\p \End(A)] = \pi_{*} (\widetilde{\alpha}^* \ch(\mathcal{O}_{\p \Sym^d(A^\vee)}(1)) \cap [\widetilde{S}]),
            \end{equation*}
            where $\widetilde{\alpha}: \widetilde{S} \to \p \Sym^d(A^\vee)$ is the regular map such that $\alpha = \pi \circ \widetilde{\alpha}$.
        \end{remark}

        In \cite{aluffiFaberAPP3}, Aluffi and Faber focus on the computation of this invariant instead of the predegree polynomial since it has better structure when considering unions of curves.

    \subsection{Segre Classes}\label{sec:segre_classes}

        A central object in the study of intersection theory is the Segre class of a proper subscheme of an ambient variety. It is defined in \S~4.2 of \cite{fulton} and plays an important role in the definition of the intersection product of varieties.
        
        Let $s(B,\p\End(A))$ be the Segre class in $\p\End(A)$ of the indeterminacy scheme $B$ of the rational map $\alpha$ determined by $X$, see \S\ref{sec:point_conditions}. The coefficients of the predegree polynomial of $X$ can be expressed directly in terms of the degrees of the components of this Segre class.
         Define:
        \begin{equation}\label{eqn:pp_def4}
            a'_i = \int_{\p \End(A)} \frac{h_1^{n^2+2n-i}}{1-dh_1} \cap ([\p \End(A)]-\ell_* s(B,\p \End(A)) \otimes \mathcal{O}_{\p \End(A)}(-d))
        \end{equation}
        where $\otimes$ refers to the tensoring operation defined in \cite{aluffiTensorClasses}.
        
        \begin{proposition}\label{prop:class_graph}
            Recall the definition of $a_i$ from \eqref{eqn:pp_def1_p1}. With the defnitions as above, $a_i = a'_i$.         
        \end{proposition}
        \begin{proof}
            To simplify the notation, for an arbitrary $r \in \mathbb{Z}$, we write $\mathcal{O}(r)$ for $\mathcal{O}_{\p \End(A)}(r)$. Let $G = \sum_{i=0}^{n} a_i h_1^i \cap [\p \End(A)]$. From Proposition 3.1 in \cite{aluffiCharClass}, 
            $$\ell_*s(B,\p \End(A)) = [\p \End(A)] - c(\mathcal{O}(d))^{-1} \cap (G \otimes \mathcal{O}(d)).$$
            Rearranging, the expression becomes 
            $$G \otimes \mathcal{O}(d) = c(\mathcal{O}(d)) \cap ([\p \End(A)] - \ell_*s(B,\p^n)).$$ 
            Tensoring both sides by $\mathcal{O}(-d)$ and using Proposition 2 in \cite{aluffiTensorClasses}, we obtain
            $$G = (c(\mathcal{O}(d)) \cap ([\p \End(A)] - \ell_*s(B,\p \End(A)))) \otimes \mathcal{O}(-d).$$
            By proposition 1 in \cite{aluffiTensorClasses},
            \begin{multline*}
                (c(\mathcal{O}(d)) \cap ([\p \End(A)] - \ell_*s(B,\p^n))) \otimes \mathcal{O}(-d) = \\ c(\mathcal{O}(-d))^{-1} \cap ([\p \End(A)] - \ell_*s(B,\p^n)) \otimes \mathcal{O}(-d)).
            \end{multline*}
            Thus,
            \begin{align*}
                a_i =& \int_{\p \End(A)} h_1^{n^2+2n-i} \cap G \\
                =& \int_{\p \End(A)}h_1^{n^2+2n-i}c(\mathcal{O}(-d))^{-1} \cap ([\p \End(A)] - \ell_*s(B,\p^n)) \otimes \mathcal{O}(-d)) \\
                =& \int_{\p \End(A)} \frac{h_1^{n^2+2n-i}}{1-dh_1} \cap ([\p \End(A)]-\ell_* s(B,\p \End(A)) \otimes \mathcal{O}_{\p \End(A)}(-d)) \\
                = & a'_i,
            \end{align*}
            as desired.
        \end{proof}

        \begin{corollary}\label{cor:explicit_class_graph}
            Let $(s^{(0)} + s^{(1)}h_1 + \cdots + s^{(n^2+2n)} h_1^{n^2+2n}) \cap [\p \End(A)]$ be the pushforward to the Chow group of $\p \End(A)$ of the Segre class of $B$. The coefficient of $t^i$ in the predegree polynomial of $X$ is independent of $s^{(i+1)}, \dots, s^{(n^2+2n)}$.
        \end{corollary}
        
        \begin{proof}
            The statement to be proven is equivalent to the following claim:
            \begin{claim}
                If $S \in A_*(\p \End(A))$ equals $\ell_*s(B,\p \End(A))$ up to codimension $i$, then
                \begin{equation*}
                    a'_i = \int_{\p \End(A)} \frac{h_1^{n^2+2n-i}}{1-dh_1} \cap ([\p \End(A)]- S \otimes \mathcal{O}_{\p \End(A)}(-d)).
                \end{equation*}
            \end{claim}
            \textit{Proof of claim.} Since this is a degree computation, it is enough to show that the difference between the rational equivalence classes whose degrees are being taken lives in lower dimension than zero; that is, that the codimension of the class
            \begin{multline}\label{eq:lowdimclass}
                \frac{1}{1-dh_1} \cap ([\p \End(A)]-\ell_* s(B,\p \End(A)) \otimes \mathcal{O}_{\p \End(A)}(-d)) \\ - \frac{1}{1-dh_1} \cap ([\p \End(A)]- S \otimes \mathcal{O}_{\p \End(A)}(-d))
            \end{multline}
            is greater than $i$. Since $A_*(\p \End(A))$ is a $A^*(\p \End(A))$-module and tensoring is linear, (\ref{eq:lowdimclass}) reduces to
            \begin{equation}\label{eq:lowdimdiff}
                \frac{1}{1-dh_1} \cap (S - \ell_*s(B,\p \End(A))) \otimes \mathcal{O}_{\p \End(A)}(-d).
            \end{equation}
            Since $S$ agrees with $\ell_*s(B,\p\End(A))$ up to codimension $i$, their difference lives in codimension greater than $i$. Then, since tensoring and multiplying by $(1-dh_1)^{-1}$ do not decrease codimension, the codimension of (\ref{eq:lowdimdiff}) is also greater than $i$ and the claim is proven.
        \end{proof}

        \begin{remark}\label{rem:independent_codim}
            Since $a_i$ and $a'_i$ can be defined at the generality of an arbitrary rational map $\rho:\p^n \dashrightarrow \p^m$, these results about the multidegrees of the graph of a rational map hold at that level of generality as well. In particular, the $i$th multidegree of the closure $\overline{\Gamma}$ of the graph of a rational map $\rho:\p^n \dashrightarrow \p^m$ is independent of the components of codimension greater than $i$ of the Segre class of its indeterminacy scheme.
        \end{remark}
        
        \begin{remark}\label{rem:small_modifications}
            Knowing components of codimension $\leq i$ of the Segre class is enough to compute the coefficient of $t^i$ in the predegree polynomial. This allows for small enough modifications to the indeterminacy scheme, which can make the computation of the coefficients of the predegree polynomial significantly simpler. We make this rigorous in the remainder of this section.
        \end{remark}
        
        Let $C$ be a closed subscheme of $\p \End(A)$. Let $U = \p \End(A) \setminus C$,
        $\ell'$ be the inclusion of $B \cap U$ in $U$ and $j'$ be the inclusion from $U$ to $\p \End(A)$. Let $g'$ denote the restriction of $j'$ to $B \cap U$.
        \begin{center}
            \begin{tikzcd}
                U \arrow[rr, "j'"] & & \p \End(A) \\
                B \cap U \arrow[rr, "g'"] \arrow[u, "\ell'"] &  & B \arrow[u, "\ell"]
            \end{tikzcd}
        \end{center}
        By Proposition 1.8 \S 1.8 in \cite{fulton}, for all $l \geq 0$ there is an exact sequence
        \begin{equation}\label{eqn:excission}
            \begin{tikzcd}
                A_l(B \cap C) \arrow[r] & A_l(B) \arrow[r,"g'^*"] & A_l(B \cap U) \arrow[r] & 0.
            \end{tikzcd}
        \end{equation}
        Since $g'^*$ is surjective, there exists $S \in A_*(B)$ such that $g'^*S = s(B \cap U, U)$. For an arbitrary choice of such $S$ define
        \begin{equation}\label{eqn:pp_def5}
            e_{i,S} = \int_{\p \End(A)} \frac{h_1^{n^2+2n-i}}{1-dh_1} \cap ([\p \End(A)]- \ell_*S \otimes \mathcal{O}_{\p \End(A)}(-d)),
        \end{equation}
        for $i$ from $0$ to $n^2+2n$. In fact, by Proposition 4.2 (b), \S 4.2 in \cite{fulton}, $s(B,\p \End(A))$ is one such class, so $a'_i = e_{i,s(B,\p\End(A))}$.

        \begin{lemma}\label{lemma:partial_resolution}
            If $i$ is less than the codimension of every component of $C$, then $a'_i = e_{i,S}$.
        \end{lemma}
        \begin{proof}
            By definition, $g'$ is an open immersion and thus flat. By Proposition 4.2 (b), \S 4.2 in \cite{fulton}, $g'^*s(B, \p \End(A)) = s(B \cap U, U)$ and by assumption, $g'^*S = s(B \cap U, U)$, so $s(B, \p \End(A)) - S$ is in the kernel of $g'^*$. Thus, by \eqref{eqn:excission}, $s(B, \p \End(A)) - S$ comes from a class in $A_*(B \cap C)$. Therefore, the dimensions of its non-zero components are at most the maximal dimension of the components of $B \cap C$; that is, $\ell_*S$ and $\ell_*s(B,\p \End(A))$ agree up to the minimal codimension of the components of $C$. Then, by the claim in Corollary \ref{cor:explicit_class_graph}, $a'_i = e_{i,S}$.
        \end{proof}

        \begin{remark}\label{rem:bezout}
            If $\codim_{\p\End(A)} B > i$, then $a_i = d^i$.
        \end{remark}
        \begin{proof}
            The indeterminacy scheme $B$ is closed, thus, consider Lemma \ref{lemma:partial_resolution} in the case $C = B$. Then, $U = \p \End(A) \setminus B$ and $g'^*s(B,\p \End(A)) = 0$. Then, $a_i = e_{i,0}$.
                    \end{proof}

        \begin{theorem}\label{thm:partial_resolution}
            Let $C$ be a closed subscheme of $\p\End(A)$. If the minimal codimension of a component of $C$ is greater than $\dim \orb(X)$, then
            \begin{equation*}
                a_i = 
                \begin{cases}
                    e_{i,S} & \text{if } 0 \leq i \leq \dim \orb(X) \\
                    0 & \text{if } i > \dim \orb(X)
                \end{cases}
            \end{equation*}
        \end{theorem}
        \begin{proof}
            If $0 \leq i \leq \dim \orb(X)$, then $i$ is less than the minimal codimension of a component of $C$ and by Lemma \ref{lemma:partial_resolution}, $a_i = e_{i,S}$. If $i > \dim \orb(X)$, then by equation \eqref{eqn:pp_def1_p2}, $a_i = 0$.
        \end{proof}

\section{The Predegree Polynomial of a Quadric in $\p^3$} \label{sec:pp_quadric_general}

    \subsection{The $n$ Dimensional Case} \label{sec:pp_quadric_n}
        In this section, we compute the leading coefficient of the predegree polynomial of a smooth quadric hypersurface in $\p^n$. This means that $\dim \p A = n$, $d = 2$ and $\dim \p \Sym^d(A^\vee) = \frac{(n+2)(n+1)}{2} - 1$. We choose homogeneous coordinates $x_0,\dots,x_n$ on $\p A$.
        
        Up to a choice of coordinates, a smooth quadric is given by the equation $x_0^2+\cdots+x_n^2=0$, the zero locus of the standard quadratic form on $A$. It follows that the linear orbit of a smooth quadric $Q$ is the set of smooth quadrics, and therefore the closure of its orbit is $\p \Sym^2(A^\vee)$; that is, of degree 1. Furthermore, this implies that the stabilizer of a smooth quadric is (up to this choice of coordinates) $\PO(A)$, the image in $\p \GL(A)$ of the group of orthogonal transformations of $A$ naturally embeded in $\GL(A)$. There is a natural embedding of $\p \GL(A)$ in $\p \End(A)$ and a natural embedding of $\GL(A)$ in $\End(A) \subset \p(\End(A) \oplus \field)$. For the rest of this discussion, the \textit{degree} of a subgroup of $\p \GL(A)$ refers to the degree of its closure under this natural embedding in $\p \End(A)$ and the degree of a subgroup of $\GL(A)$ refers to the degree of its closure in $\p(\End(A) \oplus \field)$. The convention for subgroups of $\GL(A)$ matches the one used in \cite{degreeSO}.
    
        \begin{proposition}\label{prop:lead_coeff_pp_quadric}
            The leading coefficient of the predegree polynomial of a smooth quadric hypersurface $Q \subset \p^n$ is the degree of $\PO(A)$.
        \end{proposition}
        \begin{proof}
            By Corollary \ref{cor:interpretation}, the leading coefficient of the predegree polynomial of a quadric hypersurface $Q$ is the product of the degree of the closure of its orbit and the degree of the closure of its stabilizer. By the discussion preceding this paragraph, the closure of its orbit is of degree $1$. Since its stabilizer is $\PO(A)$, the degree of the closure of its stabilizer is the degree of the closure of $\PO(A)$ in $\p \End(A)$.
        \end{proof}
        
        The degree of $\SO(A)$, the group of orientation preserving orthogonal transformations of $A$, is
        \begin{equation} \label{eq:degreeO}
            2^{n} \det\left( \binom{2n + 2 - 2i - 2j}{n + 1 - 2i}\right)_{1 \leq i,j \leq \lfloor \frac{n+1}{2} \rfloor} ,
        \end{equation}
        by Theorem 1.1 in \cite{degreeSO}. 
    
        \begin{proposition}\label{prop:degreePO}
            The degree of $\PO(A)$ equals the degree of $\SO(A)$.
        \end{proposition}
    
        \begin{proof}
            Let $I$ denote the identity matrix in $\GL(A)$. The center of $\O(A)$ consists of the matrices in $\O(A)$ of the form $aI$. Then, since $aI \in \O(A)$,
            \begin{equation*}
                I = (aI)^t aI = (aI)^2 = a^2 I.
            \end{equation*}
            Thus, $a^2 = 1$ and, since the characteristic of $\field$ is zero, $a \in \{-1,1\}$. Thus, the projection from $\O(A)$ to $\PO(A)$ is $2$ to $1$ and the projective embeddings are compatible, thus $2 \deg(\overline{\PO(A)}) = \deg(\overline{\O(A)})$. Finally, $\SO(A)$ is an index $2$ subgroup of $\O(A)$, so $2 \deg(\overline{\SO(A)}) = \deg(\overline{\O(A)})$.
        \end{proof}
    
        \begin{proposition}\label{prop:maximal_components}
            The maximal dimension of a component of the indeterminacy scheme of a general smooth quadric $Q$ in $\p^n$ is 
            \begin{equation} \label{eq:maximal_dimension_component}
                \left(n - 1 - \frac{3}{2}\left\lfloor \frac{n-1}{2} \right\rfloor\right)\left(\left\lfloor \frac{n-1}{2} \right\rfloor + 1\right) + (n + 1) \left(\left\lfloor \frac{n-1}{2} \right\rfloor + 1\right) - 1.
            \end{equation}
        \end{proposition}
        \begin{proof}
            By Theorem 22.13 in \cite{harris}, the dimension of the Fano scheme of $k$ dimensional linear subspaces of a nonsingular $n-1$-dimensional quadric, $F_{k}(Q)$ is 
            \begin{equation*}
                \left(n - 1 - \frac{3}{2}k\right)\left(k + 1\right),
            \end{equation*}
            for $k \leq \frac{n-1}{2}$, and is empty otherwise. Furthermore, if $k < \frac{n-1}{2}$, then $F_k(Q)$ is irreducible and, by Theorem 22.14 in \cite{harris}, if $n$ is odd, then it consists of two irreducible connected components. Note that, as a function of $k$, $\dim F_k(Q)$ is increasing; thus, we focus on the case $k = \left\lfloor\frac{n-1}{2} \right\rfloor$.
    
            A dimension count shows that the space of projective matrices in $\p \End(A)$ mapping to a given $k$-dimensional linear subspace of $\p A$ is $(n+1)(k+1) - 1$ dimensional. Note that, as a function of $k$, this is also increasing.
    
            Denote by $p_1$ (resp. $p_2$) the projection of $F_k(Q) \times \p \End(A)$ onto its first (resp. second) factor. Consider the subscheme $\mathcal{F}_{k}$ of $F_{k}(Q) \times \p \End(A)$ consisting of pairs $(L,\phi)$ where $\im \phi \subset L$. By Proposition \ref{prop:base_scheme}, the indeterminacy scheme $B$ of $\alpha$ is supported on the space of projective matrices whose image is contained in $Q$, so the components of $B_0$, the support of $B$, correspond to the image of $\mathcal{F}_{k}$ under the projection $p_2$.
    
            Consider the subset $U_k$ of $\mathcal{F}_{k}$ containing the projective matrices of rank $k$ whose image is in $F_k(Q)$; that is, $U_k = p_2^{-1}(R_k) \cap \mathcal{F}_k$, where $R_k$ is the open subscheme of $\p \End(A)$ of projective matrices of rank at least $k$. Then, $U_k$ is open in $\mathcal{F}_k$ and $\dim \mathcal{F}_k = \dim U_k$.
            
            Assume there exist two $k$-dimensional linear subspaces $L$ and $L'$ of $Q$ such that $(L,\phi), (L',\phi) \in U_k$ for some $\phi \in \p \End(A)$. By definition, this means that $L$ and $L'$ are both the image of $\phi$, and thus, equal. Therefore, the projection $U_k \to \p \End(A)$ is injective so that $\dim U_k = \dim p_2(U_k)$ and $\dim p_2(\mathcal{F}_k) = \dim U_k$.
    
            \begin{center}
                \begin{tikzcd}
                    & & F_k(Q) \times \p \End(A) \arrow[ddll, "p_1"'] \arrow[ddrr, "p_2"] & & \\
                    & & \mathcal{F}_{k} \arrow[u,hook] \arrow[dll] \arrow[drr] & & \\
                    F_k(Q) & & U_{k} \arrow[u,hook] \arrow[rr,hook] \arrow[ll,two heads] & & \p \End(A) \\
                \end{tikzcd}
            \end{center}

            The fiber $p_1^{-1}(L) \cap U_k$ over a given $k$-dimensional linear subspace $L \in F_k(Q)$ is always of dimension $(n+1)(k+1) - 1$ and a local trivialization exists, so $p_1:U_k \to F_k(Q)$ is a fiber bundle with constant dimensional fibers. Thus, the dimension of $U_{k}$ is the dimension of $F_k(Q)$ plus the dimension of the space of maximal rank matrices whose image is a fixed $L \in F_k(Q)$. Thus:
            \begin{align*}
                \dim p_2(\mathcal{F}_{k}) = & \dim U_k \\
                = & \dim F_k(Q) + \dim \{\phi \in \p \End(A) | \im \phi = L\} \\
                = & \left(n - 1 - \frac{3}{2}k\right)\left(k + 1\right) + (n+1)(k+1) - 1.
            \end{align*}
            Since the maximal dimension of a linear subspace of $Q$ is $\left\lfloor \frac{n-1}{2} \right\rfloor$, the maximum value of $\dim U_{k}$ is achieved when $k = \left\lfloor \frac{n-1}{2} \right\rfloor$, in which case it is \eqref{eq:maximal_dimension_component}. Thus, the maximal dimension of a component of $B$ is the maximal dimension of $\mathcal{F}_{k}$; that is, \eqref{eq:maximal_dimension_component}.
        \end{proof}
        
        Using Remark \ref{rem:bezout}, Proposition \ref{prop:degreePO}, Proposition \ref{prop:maximal_components} and \eqref{eq:degreeO}, we get Table \ref{table1} of predegree polynomials. The missing coefficients for the case $n=3$ are computed in Theorem \ref{thm:main_text}. In the case $n = 4$, the dimension of the indeterminacy scheme is $12$ so the coefficients for $t^{12}$ and $t^{13}$ are unknown and beyond the scope of this paper.
    
        \begin{center}
            \small{
            \begin{table}
                \centering
                \begin{tabular}{|c |c | c |c |}
                    \hline
                    $n$ & $\dim \p \Sym^2(A^\vee) $& $\dim \mathcal{F}_{\left\lfloor \frac{n-1}{2} \right\rfloor}$ & Predegree of a Quadric in $\p^n$ \\ [0.5ex] 
                    \hline
                    1 & 2 & 1 & $2^0t^0 + 2^1t^1 + 2t^2$  \\ 
                    \hline
                    2 & 5 & 3 & $2^0t^0 + \cdots + 2^4t^4 + 8t^5$   \\
                    \hline
                    3 & 9 & 8 & $2^0t^0 + \cdots + 2^6t^6 + *t^7  + *t^8 + 40t^9$ \\
                    \hline 
                    4 & 14 & 12 & $2^0t^0 + \cdots + 2^{11}t^{11} + *t^{12} + *t^{13} + 384t^{14}$\\
                    \hline
                \end{tabular}
                \caption{Predegree polynomials of smooth quadrics in $\p^n$.}
                \label{table1}
            \end{table}}
        \end{center}

    \subsection{The $3$ Dimensional Case} \label{sec:pp_quadric}
        For the rest of this section, we focus on the case $n=3$. As previously mentioned, all smooth quadrics in $\p^3$ lie in the same orbit, thus we choose $Q = V(x_0x_3 - x_1x_2) \subset \p A$, the image of the standard Segre embedding $\sigma:\p^1 \times \p^1 \to \p^3$. This is the hypersurface playing the role of $X$ in the discussion from \S \ref{sec:background}. Since $Q$ is isomorphic to $\p^1 \times \p^1$, there are two disjoint families of lines lying on $Q$, usually called \textit{rulings}, those of the form $\sigma(\{p\}\times \p^1)$ and those of the form $\sigma(\p^1 \times \{p\})$, for some $p \in \p^1$. We denote by the \textit{$a$-line determined by $p$} the image under $\sigma$ of $\{p\} \times \p^1$, analogously, the \textit{$b$-line determined by $p$} is $\sigma(\p^1 \times \{p\})$. Define $Z_1 \subset \p^{15}$ to be the subscheme of matrices whose image is contained in an $a$-line and, analogously, $Z_2$ is the subscheme of matrices whose image is contained in a $b$-line. Recall the definition of the rational map $\alpha$ and indeterminacy scheme $B$ from \S\ref{sec:point_conditions}.
        
        \begin{proposition}\label{prop:supp_quadric}
            The support $B_0$ of $B$ is $Z_1 \cup Z_2$.
        \end{proposition}
        \begin{proof}
            Let $\phi \in Z_1 \cup Z_2$ be arbitrary. Then, $\im \phi \subset Q$, so by Proposition \ref{prop:base_scheme}, $\phi \in B_0$.
            
            Let $\phi \in B_0$ be arbitrary. By Proposition \ref{prop:base_scheme}, $\im \phi \subset Q$. Then, the image of $\phi$ is a linear subspace of $\p^3$ contained in $Q$. So $\im \phi$ is either a point or a line. If it is a point, then $\phi \in Z_1 \cap Z_2$. If it is a line, then it is either an $a$-line or a $b$-line; thus, $\phi \in Z_1$ or $\phi \in Z_2$.
        \end{proof}
        
        A better characterization of the components of $B_0$ is given in Proposition \ref{prop:components_base_scheme}. Let $\sigma_1 : \p^1 \times \p^7 \to \p^{15}$ be the Segre embedding given by
        \begin{equation}\label{eqn:def_sigma_1}
            ((s_0:s_1),(t_0:\dots:t_7)) \mapsto 
            \begin{pmatrix}
                s_0t_0 & s_0t_1 & s_0t_2 & s_0t_3  \\
                s_0t_4 & s_0t_5 & s_0t_6 & s_0t_7  \\
                s_1t_0 & s_1t_1 & s_1t_2 & s_1t_3  \\
                s_1t_4 & s_1t_5 & s_1t_6 & s_1t_7  
            \end{pmatrix}
        \end{equation}
        and $\sigma_2: \p^1 \times \p^7 \to \p^{15}$ be given by
        \begin{equation}\label{eqn:def_sigma_2}
            ((s_0:s_1),(t_0,\dots,t_7)) \mapsto \begin{pmatrix}
                s_0t_0 & s_0t_1 & s_0t_2 & s_0t_3  \\
                s_1t_0 & s_1t_1 & s_1t_2 & s_1t_3  \\
                s_0t_4 & s_0t_5 & s_0t_6 & s_0t_7  \\
                s_1t_4 & s_1t_5 & s_1t_6 & s_1t_7
            \end{pmatrix}.
        \end{equation}
    
        The $\p^7$ factor can be identified with the space $\p M_{2,4}(\field)$ of non-zero $2 \times 4$ matrices with entries in $\field$ up to scalar multiplication.
        \begin{proposition}\label{prop:components_base_scheme}
            The components $Z_1$ and $Z_2$ of $B_0$ are the image of $\p^1 \times \p^7$ under two Segre embeddings; namely, $\sigma_1,\sigma_2$ defined by \eqref{eqn:def_sigma_1} and \eqref{eqn:def_sigma_2}.
        \end{proposition}
        \begin{proof}
            A simple computation shows that the image of an arbitrary $\phi \in \im \sigma_1 $ is a line contained in $Q$, specifically, the $a$-line determined by $(s_0:s_1)$, so $\im \sigma_1 \subset Z_1$.
            
            Now, let $\phi \in Z_1$ be arbitrary. Then, there exists $(s_0:s_1) \in \p^1$ such that $\im \phi $ is contained in the $a$-line determined by $(s_0 : s_1)$ and, as polynomials on $x_0,\dots,x_1$:
            \begin{equation}\label{eqn:a-line_containment}
                s_1(\sum_{i=0}^3 \phi_{0,i}x_i ) - s_0(\sum_{i=0}^3 \phi_{2,i}x_i ) = s_1(\sum_{i=0}^3 \phi_{1,i}x_i ) - s_0(\sum_{i=0}^3 \phi_{3,i}x_i ) = 0
            \end{equation}
            where the $\phi_{i,j}$ are the entries of $\phi$. Equation \eqref{eqn:a-line_containment} is equivalent to a system of equations which reduces to $(s_0:s_1) = (\phi_{0,i}:\phi_{2,i}) = (\phi_{1,i}:\phi_{3,i})$, which means that there is some point $(t_0:\dots:t_7) \in \p^7$ such that $\phi$ is of the form
            \begin{equation*}
                \begin{pmatrix}
                s_0 t_0 & s_0t_1 & s_0t_2 & s_0t_3  \\
                s_0t_4 & s_0t_5 & s_0t_6 & s_0t_7  \\
                s_1t_0 & s_1t_1 & s_1t_2 & s_1t_3  \\
                s_1t_4 & s_1t_5 & s_1t_6 & s_1t_7  
            \end{pmatrix}.
            \end{equation*}
            That is, $\phi \in \im \sigma_1$ and therefore $Z_1 = \im \sigma_1$.
            The proof for $Z_2$ is the same, mutatis mutandis.
        \end{proof}
        
        \begin{remark}\label{rmk:local_Z_i}
            Let $(p,\xi) \in \p^1 \times \p^7$ and denote $ \sigma_1(p,\xi)$ by $\phi_1$ and $\sigma_2(p,\xi)$ by $\phi_2$. The image of $\phi_1$ is $\sigma(\{p\} \times \im \xi)$, the image of $\phi_2 $ is $ \sigma(\im \xi \times \{p\})$ and $\p \ker \phi_1 = \p \ker \phi_2 = \p \ker \xi$. Furthermore, $\sigma_1(\p^1 \times \{\xi\})$ consists of the matrices $\zeta$ such that there exists $r \in \p^1$ for which $\im \zeta = \sigma(\{r\} \times \im \xi)$ and $\p \ker \xi \subset \p \ker \zeta$. The image of $\{p\} \times \p^7$ under $\sigma_1$ consists of the matrices whose image is contained in the $a$-line determined by $p$. The description of $\sigma_2(\p^1 \times \{\xi\})$ and $\sigma_2(\{p\} \times \p^7)$ is analogous. Finally, $\sigma_1(p,\xi),\sigma_2(p,\xi) \in Z_1 \cap Z_2$ if and only if the rank of $\xi$ is $1$.
        \end{remark}
    
        The inverse of the Chern class of the normal bundle $N$ to the image of a Segre embedding $\sigma': \p^{n_1} \times \dots \times \p^{n_r} \to \p^m$ is
        \begin{equation}\label{eqn:segre_segre}
            c(N)^{-1} = \frac{\Pi^{r}_{i=1} (1+h_i)^{n_i+1}}{(1 + \sum_{i=1}^r h_i)^{m+1}},
        \end{equation}
        where $h_i$ is (the pullback by the $i$th projection of) $c_1(\mathcal{O}_{\p^{n_i}}(1))$ (cf.~\cite[Example~3.2.15(c)]{fulton}). Using \eqref{eqn:segre_segre} and the fact that, by Proposition 4.1 \S 4.1 in \cite{fulton}
        \begin{equation}\label{eqn:segre_from_chern}
            s(X,\p^m) = c(N)^{-1} \cap [\im \sigma'],
        \end{equation}
        we compute the Segre classes of $Z_1$ and $Z_2$.
    
        \begin{proposition}\label{prop:class_Z_i}
            The pushforward by inclusion of the Segre class of the image $Z \cong \p^1 \times \p^7$ of a Segre embedding in $\p^{15}$ is $(8H^7 - 70H^8 + 344H^9 - 1248H^{10} + 3720H^{11} - 9636H^{12} + 22440H^{13} - 48048H^{14} + 96096H^{15}) \cap [\p^{15}]$. In particular, this equals the push-forward of $s(Z_1,\p^{15})$ and $s(Z_2,\p^{15})$.
        \end{proposition}
        \begin{proof}
            Let $\text{pr}_1$ and $\text{pr}_2$ be the projections of $Z$ onto $\p^1$ and $\p^7$, respectively; $j$ be the inclusion of $Z$ in $\p^{15}$ and $N$ denote the normal bundle to $Z$ in $\p^{15}$. Also, let $H = c_1(\mathcal{O}_{\p^{15}}(1)), k_1 = \text{pr}_1^* c_1(\mathcal{O}_{\p^1}(1))$, $k_2 = \text{pr}_2^* c_1(\mathcal{O}_{\p^7}(1))$. Then, \eqref{eqn:segre_segre} reduces to:
            \begin{equation*}
                c(N)^{-1} = \frac{(1+k_1)^2(1+k_2)^8}{(1+k_1+k_2)^{16}}
            \end{equation*}
            and \eqref{eqn:segre_from_chern} becomes
            \begin{multline*}
                s(Z,\p^{15}) = c(N)^{-1} \cap [Z] = ( 1 - 14k_1 - 8k_2 + 128k_1k_2 + 36k_2^2 - 648k_1k_2^2 - 120k_2^3 \\+ 2400k_1k_2^3 + 330k_2^4 - 7260k_1k_2^4 - 792k_2^5 + 19008k_1k_2^5 + 1716k_2^6 - 44616k_1k_2^6 \\ - 3432k_2^7 + 96096k_1k_2^7 ) \cap [Z].
            \end{multline*}
            Note that $j^* H = k_1 + k_2$. Therefore, for all pairs of non-negative integers $m,n$:
            \begin{equation*}
                j_* (k_1^n k_2^m \cap [Z]) = \binom{8 - n - m}{7 - m}H^{7+n+m} \cap [\p^{15}].
            \end{equation*}
            Thus:
            \begin{multline*}
                j_*s(Z,\p^{15}) = (8H^7 - 70H^8 + 344H^9 - 1248H^{10} + 3720H^{11} - 9636H^{12} \\+ 22440H^{13} - 48048H^{14} + 96096H^{15}) \cap [\p^{15}],
            \end{multline*}
            as stated.
        \end{proof}
        \begin{proposition}\label{prop:intersection_components_base_scheme}
            The support of $Z_1 \cap Z_2$ is isomorphic to $Q \times \p^3$. Therefore, the codimension of $Z_1 \cap Z_2$ in $\p^{15}$ is $10$.
        \end{proposition}
        \begin{proof}
            The locus of rank $1$ matrices in $\p^{15} = \p \End(\A^4)$ is given by a Segre embedding $\rho:\p^3 \times \p^3 \to \p^{15}$, where the first factor parameterizes the image in $\p^3$ of the matrix and the second factor parameterizes the dual of the projectivization of the kernel of a representative of the matrix.
    
            An element $\phi \in \p^{15}$ is in the intersection of $Z_1$ and $Z_2$ if and only if the image of $\phi$ is contained in the $a$-line and $b$-line determined by $p$ and $q$ in $\p^3$, respectively. Then, $\phi$ is in the intersection of $Z_1$ and $Z_2$ if and only if the image of $\phi$ is $\sigma(p,q)$. Thus $\phi$ is in the intersection of $Z_1$ and $Z_2$ if and only if $\phi$ is in the image of $Q \times \p^3$ under $\rho$.
        \end{proof}
    
        \begin{proposition} \label{prop:base_scheme_open_reduced}
            The indeterminacy scheme of $\alpha$ is reduced away from $Z_1 \cap Z_2$.
        \end{proposition}
    
        \begin{proof}
            Consider $\phi \in Z_1 \setminus Z_2$. By Proposition \ref{prop:components_base_scheme}, $\phi$ corresponds to a point $(p,\xi) \in \p^1 \times \p^7$ and since $\phi \notin Z_1 \cap Z_2$, the rank of $\xi$ is $2$. Then, up to a choice of coordinates, $p = (1:0)$ and 
            \begin{equation*}
                \xi = 
                \begin{pmatrix}
                    1 & 0 & 0 & 0 \\
                    0 & 1 & 0 & 0
                \end{pmatrix}.
            \end{equation*}
            Thus, $\p \ker \phi$ is $V(x_0,x_1) \subset \p^3$ and for $q = (q_0:q_1:q_2:q_3) \notin \p \ker \phi$, $\phi(q) = (q_0:q_1:0:0)$.
            
            Let $q = (q_0:q_1:q_2:q_3) \in \p^3 \setminus \p \ker \phi$ be arbitrary. Then, by Proposition \ref{prop:tangent_point_condition}, the embedded tangent space to $P_q$ at $\phi$ is
            \begin{equation*}
                \{\psi \in \p^{15} | q \in \p \ker \psi \text{ or } \psi(q) \in \TEmb_{\phi(q)} Q\}.
            \end{equation*}
            Then, the ideal corresponding to $\TEmb_\phi P_q$ generated by
            \begin{equation*}
                \tau_q = q_0 \sum_{i=0}^3 a_{3,i} q_i - q_1 \sum_{i=0}^3a_{2,i} q_i,
            \end{equation*}
            where the $a_{i,j}$ are the homogeneous coordinates on $\p^{15}$. Recall that $B$, the indeterminacy scheme of $\alpha$ is the (scheme-theoretic) intersection of all point conditions. Then, 
            \begin{equation*}
                \TEmb_{\phi} B = \bigcap_{q \in \p^3} \TEmb_\phi P_q,
            \end{equation*}
            so $\TEmb_\phi B = V(\{\tau_q| q \in \p^3\})$. The ideal generated by $\{\tau_q| q \in \p^3\}$ is also generated by 
            \begin{equation*}
                \{a_{2,0} - a_{3,1}, a_{2,1},a_{2,2},a_{2,3},a_{3,0},a_{3,2},a_{3,3}\}.
            \end{equation*}
            Therefore, $\TEmb_\phi B$ is $8$ dimensional and it contains $\TEmb_\phi Z_1$ which is also $8$ dimensional. Thus, $\TEmb_\phi B = \TEmb_\phi Z_1$ and $B$ is reduced over $Z_1 \setminus Z_2$. Analogously, $B$ is reduced over $Z_2 \setminus Z_1$.
        \end{proof}
        
        Let $U$ denote $\p^{15} \setminus (Z_1 \cap Z_2)$. Let $\ell'$ the inclusion of $B \cap U$ in $U$, $j'$ the inclusion of $U$ in $\p^{15}$, $i_1$, resp. $i_2$, the inclusion of $Z_1$, resp. $Z_2$ in $B$. Let $g'$ denote the restriction of $j'$ to $B \cap U$. This gives the following commutative diagram.
        \begin{center}
            \begin{tikzcd}
                & & & \p^{15} & \\
                & U \arrow[rru, "j'"] & & & \\
                & & Z_2 \arrow[r, "i_2"] \arrow[ruu, "\ell_2"] & B \arrow[uu, "\ell"] & Z_1 \arrow[l, "i_1"'] \arrow[luu, "\ell_1"'] \\
                Z_2 \cap U \arrow[r, "i'_2"'] \arrow[rru, "g'_2"'] \arrow[ruu, "\ell'_2"'] & B \cap U \arrow[rru, "g'"'] \arrow[uu, "\ell'"'] & Z_1 \cap U \arrow[l, "i'_1"] \arrow[rru, "g'_1"'] \arrow[luu, "\ell'_1"'] & &
            \end{tikzcd}
        \end{center}
        \begin{proposition}\label{prop:s_works}
            Let $S = i_{1*}s(Z_1,\p^{15}) + i_{2*}s(Z_2,\p^{15})$. The pullback by $g'$ of $S$ is $s(B \cap U,U)$.
        \end{proposition}
        \begin{proof}
            Since $B\cap U = Z_1 \cap U \sqcup Z_2 \cap U$, by Example 1.3.1(b) \S 1.3 in \cite{fulton}, $A_*(B \cap U) \cong A_*(Z_1 \cap U) \oplus A_*(Z_2 \cap U)$, where the isomorphism is given by ${i'_1}^* \oplus {i'_2}^*$. By Proposition 4.2(b) \S 4.2 in \cite{fulton}, ${i'_1}^*s(B \cap U,U) = s(Z_1 \cap U,U)$ and ${i'_2}^*s(B \cap U,U) = s(Z_2 \cap U,U)$, thus
            \begin{align*}
                ({i'_1}^* \oplus {i'_2}^*)(s(B\cap U,U)) = & ({i'_1}^*s(B\cap U,U), {i'_2}^*s(B\cap U,U))\\
                = & (s(Z_1 \cap U,U) , s(Z_2 \cap U,U)).
            \end{align*}
            Also by Proposition 4.2(b) \S 4.2 in \cite{fulton},
            \begin{align*}
                {g'}_1^*s(Z_1,\p^{15}) = & s(Z_1 \cap U,U) \\
                {g'}_2^*s(Z_2,\p^{15}) = & s(Z_2 \cap U,U).
            \end{align*}
            Thus, $g'^*S = (s(Z_1 \cap U,U),s(Z_2 \cap U,U)) = s(B\cap U,U)$.
        \end{proof}
    
        \begin{corollary}\label{cor:partial_segre}
            The pushforward by $\ell$ of $S$ is $2\ell_{1*}s(Z_1,\p^{15})$.
        \end{corollary}
        \begin{proof}
            Since $\ell \circ i_1 = \ell_1$ and $\ell \circ i_2 = \ell_2$, and, by Proposition \ref{prop:class_Z_i}, $\ell_{2*} s(Z_2,\p^{15}) = \ell_{1*}s(Z_1,\p^{15})$,
            \begin{align*}
                \ell_* S = & \ell_*i_{1*} s(Z_1,\p^{15}) + \ell_* i_{2*} s(Z_2,\p^{15}) \\
                = & \ell_{1*} s(Z_1,\p^{15}) + \ell_{2*}s(Z_2,\p^{15}) \\
                = & \ell_{1*} s(Z_1,\p^{15}) + \ell_{1*} s(Z_1,\p^{15}) \\
                = & 2\ell_{1*} s(Z_1,\p^{15}),
            \end{align*}
            as stated.
        \end{proof}
        
        \begin{theorem}\label{thm:main_text}
            The predegree polynomial of a smooth quadric $Q$ in $\p^3$ is \hfill \break 
            $\mathcal{P}_Q(t) = 1+ 2t + 4t^2 + 8t^3 + 16t^4 + 32t^5 + 64t^6 + 112t^7 + 140t^8 + 40t^9$
        \end{theorem}
        \begin{proof}
            By Proposition \ref{prop:intersection_components_base_scheme}, the codimension of $Z_1 \cap Z_2$ in $\p^{15}$ is larger than $\dim \orb(Q) = 9$. Then, by Theorem \ref{thm:partial_resolution}, for $i$ from $0$ to $9$ $a_i = e_{i,S}$, where $S$ is the class defined in Proposition \ref{prop:s_works}. The computation is a straightforward application of Corollary \ref{cor:partial_segre} and \eqref{eqn:pp_def5}.
        \end{proof}
    
        Theorem \ref{thm:main_text} completes the entry for $n=3$ in Table \ref{table1}. Given a general $n$ dimensional linear subspace $L \subset \p^{15}$, denote by $Q \circ L$ the set of quadrics $Q'$ such that for some invertible $\phi \in L$, $Q'= Q \circ \phi$. Using this notation, Corollary \ref{cor:interpretation} interprets Theorem \ref{thm:main_text}  as statements of the form `given a general $7$ dimensional linear subspace $L$ of $\p^{15}$ and $7$ points in general position in $\p^3$, there are 112 quadrics $Q'$ in $Q \circ L$ such that $Q'$ contains all $7$ points.' Note that the case for $\dim L = 9$ is special since the translates here are counted with multiplicity $\deg \PO(4)$. These statements are summarized in Table \ref{table:interpretation}.
        \begin{center}
            \begin{table}
                \centering
                \small{
                \begin{tabular}{|c |c | c |c |} 
                    \hline
                    $\dim L =$ \# of points & \# of $Q' \in Q \circ L$ containing all points \\ [0.5ex] 
                    \hline
                    0 & 1 \\
                    \hline
                    1 & 2 \\ 
                    \hline
                    2 & 4 \\
                    \hline
                    3 & 8 \\
                    \hline 
                    4 & 16 \\
                    \hline
                    5 & 32 \\
                    \hline
                    6 & 64 \\
                    \hline
                    7 & 112 \\
                    \hline
                    8 & 140 \\
                    \hline
                    9 & 1 \\
                    \hline
                \end{tabular}}
                \caption{Interpretation of the main result}
                \label{table:interpretation}
            \end{table}
        \end{center}

    \section{Resolving the Rational Map} \label{sec:resolution}

    In this section, we resolve the indeterminacies of the map $\alpha$ from \S \ref{sec:point_conditions} in the specific case where $A$ is a $4$-dimensional vector space and $X$ is the smooth quadric $Q$ defined in \S\ref{sec:pp_quadric}. As proven in Proposition \ref{prop:supp_quadric}, the support of the indeterminacy scheme of $\alpha$ is $Z_1 \cup Z_2$.

    \subsection{Tangent Spaces}

         \begin{remark} \label{rmk:tangent_Z_i}
            The embedded tangent space to $Z_1$ at $\phi = \sigma_1(p,\xi)$ is the $\p^8$ spanned by $\sigma_1(\{p\} \times \p^7)$ and $\sigma_1(\p^1 \times \{\xi\})$. Analogously, The embedded tangent space to $Z_2$ at $\phi = \sigma_2(p,\xi)$ is the $\p^8$ spanned by $\sigma_2(\{p\}\times \p^7)$ and $\sigma_2(\p^1 \times \{\xi\})$. These spaces are described in Remark \ref{rmk:local_Z_i}.
        \end{remark}
    
        Consider $\phi \in Z_1 \cap Z_2$. Then, by Proposition~\ref{prop:components_base_scheme}, there is a unique point $(p,q,k) \in \p^1 \times \p^1 \times \p^3 \cong Q \times \p^3$ such that $\im \phi \subset \sigma(p,q)$ and $\p \ker \phi = k^\vee$, the plane dual to $k$. Let $L_{p,a}$ denote the $a$-line determined by $p$ and $L_{q,b}$ denote the $b$-line determined by $q$. 
    
        Consider the standard Segre embedding $\rho':\p^1 \times \p^3 \to \p^7$. A quick check shows that for all $(p,q,k) \in \p^1 \times \p^1 \times \p^3$,
        \begin{equation*}\label{eqn:funceqn_1}
            \rho(\sigma(p,q),k) = \sigma_1(p,\rho'(q,k))
        \end{equation*}
        \begin{equation*}\label{eqn:funceqn_2}
            \rho(\sigma(p,q),k)= \sigma_2(q,\rho'(p,k))
        \end{equation*}
        where $\rho$ was defined in the proof of Proposition \ref{prop:intersection_components_base_scheme}, $\sigma$ is the standard Segre embedding of $\p^1 \times \p^1$ in $\p^3$ and $\sigma_1$ (resp. $\sigma_2$) was defined in \eqref{eqn:def_sigma_1} (resp.~\eqref{eqn:def_sigma_2}). Under the identification $\p^7 \cong \p M_{2,4}(\field)$, $\p \ker \rho'(p,k) = k^\vee$ and $\im \rho'(p,k) = p$.
    
        \begin{proposition}\label{prop:tangent_Z_i}
            The embedded tangent space to $Z_1$ at $\phi = \rho(\sigma(p,q),k)$ is the $8$ dimensional linear subspace of $\p^{15}$ containing the matrices whose image is contained in $L_{p,a}$ and the matrices whose image is a point in $L_{q,b}$ and whose kernel is $k^\vee$. Similarly, the embedded tangent space to $Z_2$ at $\phi \in Z_1 \cap Z_2$ is the $\p^8$ containing the matrices whose image is contained in $L_{q,b}$ and matrices whose image is a point in $L_{p,a}$ and whose kernel is $k^\vee$.
        \end{proposition}
        \begin{proof}
            By Remarks \ref{rmk:local_Z_i} and \ref{rmk:tangent_Z_i}, the tangent space to $Z_1$ at $\phi = \sigma_1(p,\rho'(q,k))$ is given by the span of the matrices $\psi$ such that $\im \psi \subset L_{p,a}$ and the matrices $\zeta$ such that $\p \ker \zeta \supseteq \p \ker \rho'(q,k) = k^\vee$ and, for some $r \in \p^1$, the image of $\zeta$ is $\sigma(\{r\} \times \im \xi) = \sigma(\{r\} \times \{q\}) \subset L_{q,b}$. Note that since the zero matrix is not an element of $\p^{15}$ and $\codim_{\p^3} k^\vee = 1$, so if $\p \ker \zeta \supseteq k^\vee$, then $\p \ker \zeta = k^\vee$.
            
            The proof for the tangent space to $Z_2$ at $\phi$ is analogous.
        \end{proof}
    
        \begin{proposition}\label{prop:tangent_C_1}
            The tangent space to the smooth variety $Z_1 \cap Z_2$ at $\phi$, corresponding to a point $(p,q,k) \in \p^1 \times \p^1 \times \p^3$ is the space of matrices such that $\phi(k^\vee) = \sigma(p,q)$ or $k^\vee \subset \p \ker \phi$ and $\im \phi$ is contained in the plane spanned by $L_{p,a}$ and $L_{q,b}$.
        \end{proposition}
        \begin{proof}
            A computation similar to the one in the proof of the previous proposition shows this.
        \end{proof}
    
        \begin{corollary}\label{cor:intersection_tangents}
            For all $\phi \in Z_1 \cap Z_2$, the intersection of the tangent space to $Z_1$ at $\phi$ and the tangent space to $Z_2$ at $\phi$ equals the tangent space to $Z_1 \cap Z_2$ at $\phi$.
        \end{corollary}
        \begin{proof}
            Up to a choice of coordinates, any choice of $\phi \in Z_1 \cap Z_2$ is equivalent, so assume that $\phi = \rho((1:0),(0:1),(0:0:1:0))$. Then
            \begin{align*}
                \TEmb_\phi Z_1 = &V(a_{0,0},a_{0,1},a_{0,3},a_{1,0},a_{1,1},a_{1,2},a_{1,3}) \\
                \TEmb_\phi Z_2 = &V(a_{1,0},a_{1,1},a_{1,2},a_{1,3},a_{3,0},a_{3,1},a_{3,3}) \\
                \TEmb_\phi (Z_1 \cap Z_2) = & V(a_{0,0},a_{0,1},a_{0,3},a_{1,0},a_{1,1},a_{1,2},a_{1,3},a_{3,0},a_{3,1},a_{3,3}),
            \end{align*}
            so $\T_\phi Z_1 \cap \T_\phi Z_2 = \T_\phi (Z_1 \cap Z_2)$.
        \end{proof}
        
    \subsection{The First Blow-up}

    Let $C_1$ denote $Z_1 \cap Z_2 \cong Q \times \p^3$. Let $\pi_1: V_1 \to \p^{15}$ be the blow-up of $\p^{15}$ along $C_1$, $\widetilde{Z_1}$ be the proper transform of $Z_1$ and $\widetilde{Z_2}$ be the proper transform of $Z_2$. Let $E_1 := \pi_1^{-1} (C_1)$; that is, the exceptional divisor in $V_1$. The blow-up $V_1$ gives rise to the following diagram:
    \begin{center}
        \begin{tikzcd}
            V_1 \arrow[drr,"\alpha_1",dashed] \arrow[d,"\pi_1"] & & \\
            \p^{15} \arrow[rr,"\alpha",dashed] & & \p^9
        \end{tikzcd}
    \end{center}
    where $\pi_1$ is the blow-up map, $\alpha$ is the rational map defined in \S \ref{sec:point_conditions} and $\alpha_1$ is the rational map induced on $V_1$ by $\alpha$.
    
    \begin{proposition}\label{prop:first_blow-up}
        The intersection of $\widetilde{Z_1}$ and $\widetilde{Z_2}$ is empty.
    \end{proposition}
    \begin{proof}
        Since $C_1 = Z_1 \cap Z_2$, $\widetilde{Z_1} \cap \widetilde{Z_2} \subset E_1$. Note that $\widetilde{Z_i} \cap E_1$ is the image of the exceptional divisor in the blow-up of $Z_i$ along $C_1$ through the embedding $\widetilde{Z_i} \xhookrightarrow{} V_1$. This allows for the identification of points in $\widetilde{Z_i} \cap E_1$ with normal directions to $C_1$ that are tangent to $Z_i$; that is, we identify $\widetilde{Z_i} \cap E_1$ with the projectivization of the normal bundle to $C_1$ in $Z_i$, namely, $\N_{C_1} Z_i$. 
        
        The statement $\widetilde{Z_1} \cap \widetilde{Z_2} = \emptyset$ is equivalent to the statement that, as subbundles of $\N_{C_1} \p^{15}$, $ \N_{C_1} Z_1 \cap \N_{C_1} Z_2 = C_1 \times \{0\}$; that is, for every $\phi \in C_1$, $\T_\phi Z_1 / \T_\phi C_1  \cap \T_\phi Z_1 / \T_\phi C_1 = \{0\}$, as subspaces of $\T_\phi \p^{15} / \T_\phi C_1$.
        
        The statement at the level of $\T_\phi \p^{15}$ becomes $\T_\phi Z_1 \cap \T_\phi Z_2 \subset \T_\phi C_1$, which by Corollary~\ref{cor:intersection_tangents}, holds. Thus $\widetilde{Z_1} \cap \widetilde{Z_2} = \emptyset$. 
    \end{proof}
    
    Let $\mathcal{G}$ be the bundle over $C_1$ whose fiber over $\phi \in C_1$ is given by $\bigcap_{q \in \p^{3}} \T_\phi P_q$; that is, the intersection in $\T_\phi \p^{15}$ of the tangent spaces to each point condition.

    \begin{proposition}\label{prop:rank_G}
        The rank of $\mathcal{G}$ is $11$.
    \end{proposition}

    \begin{proof}
        Up to a choice of coordinates, any choice of $\phi$ is equivalent to 
        \begin{equation*}
            \begin{pmatrix}
                1 & 0 & 0 & 0 \\
                0 & 0 & 0 & 0 \\
                0 & 0 & 0 & 0 \\
                0 & 0 & 0 & 0 
            \end{pmatrix}.
        \end{equation*}
        Let $q \in \p^3$ be arbitrary. By Proposition \ref{prop:tangent_point_condition}, 
        \begin{equation*}
            \TEmb_\phi P_q = V(q_0 \sum_{i=0}^3 q_i a_{3,i}).
        \end{equation*}
        Thus,
        \begin{align*}
            \bigcap_{q \in \p^{3}} \TEmb_\phi P_q = & \bigcap_{q \in \p^{3}} V(q_0 \sum_{i=0}^3 q_i a_{3,i}) \\
            = & V(a_{3,0}, a_{3,1}, a_{3,2}, a_{3,3}).
        \end{align*}
        Therefore, $\rk \mathcal{G} = \dim \bigcap_{q \in \p^{3}} \T_\phi P_q = \codim_{\p^{15}} V(a_{3,0}, a_{3,1}, a_{3,2}, a_{3,3}) = 11.$
    \end{proof}
    
    The bundle $\mathcal{G}$ contains $\T C_1$ as a subbundle, thus the cokernel of the inclusion of $\T C_1$ into $\mathcal{G}$ is a vector bundle over $C_1$ which we denote $\mathcal{F}$. Let $C_2$ denote $\p (\mathcal{F})$, the projectivization of the bundle $\mathcal{F}$ over $C_1$. The projective bundle $C_2$ embeds into $\p \N_{C_1} \p^{15}$ and thus it is a subscheme of $E_1$.
    \begin{remark}\label{rem:rank_F}
        The rank of $\mathcal{F}$ is $6$ and the dimension of $C_2$ is $10$. Furthermore, $C_2$ is smooth.
    \end{remark}

    \begin{proposition}\label{lemma:supp_alpha_1}
         The support of the indeterminacy scheme of $\alpha_1$ is the union of
        $\widetilde{Z_1}$, $\widetilde{Z_2}$ and $C_2$. Furthermore, $\widetilde{Z_1} \cap C_2 = \widetilde{Z_1} \cap E_1$ and $\widetilde{Z_2} \cap C_2 = \widetilde{Z_2} \cap E_1$.
    \end{proposition}
    \begin{proof}
        It is clear that $\widetilde{Z_1}$ and $\widetilde{Z_2}$ are contained in the support of the indeterminacy scheme of $\alpha_1$. Since $\alpha_1$ is defined by the linear system consisting of the proper transforms of the point conditions, all that is left to show is that over $C_1$, the intersection of all point conditions is indeed $C_2$.

        Let $\phi \in C_1$ and $q \in \p^3 \setminus \p \ker \phi$. Consider $\widetilde{P_q} \cap E_1$, the exceptional divisor in the blow-up of $P_q$ along $C_1$. Then, in the fiber over $\phi$, $\widetilde{P_q} \cap E_1 = \p \N_{C_1} P_q$; that is, a point in $\widetilde{P_q} \cap E_1$ corresponds to a tangent direction to $P_q$ that is normal to~$C_1$. 
        
        Since the support of the indeterminacy scheme of $\alpha_1$ is the (set-theoretic) intersection of all point conditions in $V_1$, a point $\psi$ in the fiber of $E_1$ over $\phi \in C_1$ is in the support of the indeterminacy scheme if it is a direction tangent to all point conditions and normal to $C_1$ at $\phi$; that is, if its representatives are in the fiber of $\mathcal{F}$ over $\phi$ and thus $\psi$ is a point in $C_2$. Therefore the support of the indeterminacy scheme of $\alpha_1$ is $\widetilde{Z_1} \cup \widetilde{Z_2} \cup C_2$.

        Since $C_2 \subset E_1$, $\widetilde{Z_i} \cap C_2 \subset \widetilde{Z_i} \cap E_1$.
        
        Recall from the proof of Proposition \ref{prop:first_blow-up} that $\widetilde{Z_1} \cap E_1 = \p \N_{C_1} Z_1$; that is, the projectivization of the normal bundle to $C_1$ in $Z_1$. Thus, $\widetilde{Z_1} \cap E_1 \subset \widetilde{Z_1} \cap C_2$ is equivalent to the statement that, over $C_1$, $\T Z_1$ is a subbundle of $\mathcal{G}$.
        Since $Z_1$ is contained in the support of the indeterminacy scheme of $\alpha$, $\T Z_1$ is a subbundle of $\mathcal{G}$, thus $\widetilde{Z_1} \cap E_1 \subset \widetilde{Z_1} \cap C_2$.

        Analogously, $\widetilde{Z_2} \cap E_1 = \widetilde{Z_2} \cap C_2$.
    \end{proof}

    \subsection{Further Blow-ups}
     Let $\pi_2:V_2 \to V_1$ be the blow-up of $V_1$ along $C_2$. Since $\pi_2$ is birational, $\alpha_1$ induces a rational map $\alpha_2:  V_2 \dashrightarrow \p^9$. This information is captured in the following diagram:

    \begin{center}
        \begin{tikzcd}
            & V_2 \arrow[rdd, "\alpha_2", dashed] \arrow[ld, "\pi_2"] & \\
            V_1 \arrow[rrd, "\alpha_1", dashed] \arrow[d, "\pi_1"] &  &  \\
            \p^{15} \arrow[rr, "\alpha", dashed] &  & \p^{9}
        \end{tikzcd}
    \end{center}
    
    \begin{proposition}\label{prop:second_blow-up}
        The indeterminacy scheme of $\alpha_2$ is the union of $\widehat{{Z_1}}$ and $\widehat{{Z_2}}$, the proper transforms of $\widetilde{Z_1}$ and $\widetilde{Z_2}$.
    \end{proposition}
    \begin{proof}
        By Proposition \ref{prop:base_scheme_open_reduced}, away from $E_2$ and $\overline{\pi_2^{-1}(E_1 \setminus C_2)}$ the indeterminacy scheme is indeed $\widehat{Z_1} \cup \widehat{Z_2}$. Consider the affine patch in $\p \End(A) \cong \p^{15}$ where $a_{0,0}\neq 0$. We choose coordinates $(a_{0,1}, \cdots, a_{3,3})$ for $\A^{15}$, induced by the coordinates in $\p^{15}$, and $(b_0: \cdots : b_{9})$ for $\p^9$, which give coordinates to $\pi_1^{-1}(D(a_{0,0})) \subset \A^{15} \times \p^{9}$. Consider the affine patch where $b_{0} \neq 0$. In this affine patch, the change of variables formulae are:
        \begin{align*}
            a_{1,2} \mapsto & a_{0,2}a_{1,0} + b_1(a_{1,1} - a_{0,1}a_{1,0}) \\
            a_{1,3} \mapsto & a_{0,3}a_{1,0} + b_2(a_{1,1} - a_{0,1}a_{1,0}) \\
            a_{2,1} \mapsto & a_{0,1}a_{2,0} + b_3(a_{1,1} - a_{0,1}a_{1,0}) \\
            a_{2,2} \mapsto & a_{0,2}a_{2,0} + b_4(a_{1,1} - a_{0,1}a_{1,0}) \\
            a_{2,3} \mapsto & a_{0,3}a_{2,0} + b_5(a_{1,1} - a_{0,1}a_{1,0}) \\
            a_{3,0} \mapsto & a_{1,0}a_{2,0} + b_6(a_{1,1} - a_{0,1}a_{1,0}) \\
            a_{3,1} \mapsto & a_{0,1}a_{1,0}a_{2,0} + b_7(a_{1,1} - a_{0,1}a_{1,0}) \\
            a_{3,2} \mapsto & a_{0,2}a_{1,0}a_{2,0} + b_8(a_{1,1} - a_{0,1}a_{1,0}) \\
            a_{3,3} \mapsto & a_{0,3}a_{1,0}a_{2,0} + b_9(a_{1,1} - a_{0,1}a_{1,0}), \\
        \end{align*}
        and the equation of $E_1$ is $a_{1,1} - a_{0,1}a_{1,0} =0$. Using this information and the definition of $\alpha$, the indeterminacy scheme of $\alpha_1$ is determined by the ideal
        \begin{multline*}
            (b_6, b_7 - a_{1,0}b_3 - a_{2,0}, b_8 - a_{1,0}b_4 - a_{2,0}b_1, b_9 - a_{1,0}b_5 - a_{2,0}b_2, \\
            b_4(a_{1,1} - a_{0,1}a_{1,0}), b_5(a_{1,1} - a_{0,1}a_{1,0}), b_3(a_{1,1} - a_{0,1}a_{1,0})),
        \end{multline*}
        which has two components, one being the ideal of $\widetilde{Z_2}$ and the other the ideal of $C_2$. This is the case for the patches where $b_0,b_1$ or $b_2$ are non-zero. Over the patches where $b_3,b_4$ or $b_5$ are non-zero, the indeterminacy scheme has two components as well, but they are $\widetilde{Z_1}$ and $C_2$ instead. Over the patch where $b_6 \neq 0$, the indeterminacy scheme is empty and over the patches where $b_7,b_8$ or $b_9$ are non-zero, the indeterminacy scheme has three components, $\widetilde{Z_1}, \widetilde{Z_2}$ and $C_2$. This means that the indeterminacy scheme of $\alpha_2$ over $\overline{\pi_2^{-1}(E_1 \setminus C_2)}$ is $\widehat{Z_1} \cup \widehat{Z_2}$. 
        
        Using the change of variables formulae, $\pi_1^{-1}(D(a_{0,0})) \cap D(b_0) \cong \A^{15}$ with coordinates $(a_{0,1}, \cdots, a_{1,1}, a_{2,0},b_1,\cdots,b_5)$. Then, $\pi_2^{-1}(\pi_1^{-1}(D(a_{0,0})) \cap D(b_0)) \subset \A^{15} \times \p^4$ and we choose homogeneous coordinates $(c_0: \cdots:c_4)$ on $\p^4$. Consider the affine patch of $\pi_2^{-1}(\pi_1^{-1}(D(a_{0,0})) \cap D(b_0))$ where $c_4 \neq 0$. In this affine patch, the change of variables formulae are:
        \begin{align*}
            b_6 \mapsto & c_0(a_{1,1} - a_{0,1}a_{1,0}) \\ 
            b_7 \mapsto & a_{1,0}b_3 + a_{2,0} + c_1(a_{1,1} - a_{0,1}a_{1,0})  \\
            b_8 \mapsto & a_{1,0}b_4 + a_{2,0}b_1 + c_2(a_{1,1} - a_{0,1}a_{1,0}) \\ 
            b_9 \mapsto & a_{1,0}b_5 + a_{2,0}b_2 + c_3(a_{1,1} - a_{0,1}a_{1,0}),
        \end{align*}
        and the equation of $E_2$ is $a_{1,1} - a_{0,1}a_{1,0}$. Using this information and the definition of $\alpha$, the indeterminacy scheme of $\alpha_2$ is determined by the ideal
        \begin{equation*}
            (c_0,c_1,c_2,c_3,b_3,b_4,b_5),
        \end{equation*}
        which is the ideal of $\widehat{Z_2}$. Going through each patch and doing the relevant computations we obtain that the indeterminacy scheme of $\alpha_2$ is precisely $\widehat{Z_1} \cup \widehat{Z_2}$.
    \end{proof}

    Let $\pi_3:V_3 \to V_2$ be the blow-up of $V_2$ along $\widehat{{Z_1}}$ and $\widehat{{Z_2}}$. Then, since $\pi_3$ is birational, $\alpha_2$ induces a map $\alpha_3:V_3 \dashrightarrow \p^9$, which our final result shows is regular. This information is captured in the following diagram: 
    \begin{center}
        \begin{tikzcd}
            & V_2 \arrow[rdd, "\alpha_2", dashed] \arrow[ld, "\pi_2"] & V_3 \arrow[l, "\pi_3"] \arrow[dd, "\alpha_3"] \\
            V_1 \arrow[rrd, "\alpha_1", dashed] \arrow[d, "\pi_1"] &  & \\
            \p^{15} \arrow[rr, "\alpha", dashed] &  & \p^{9}
        \end{tikzcd}
    \end{center}
   
    \begin{theorem}\label{thm:third_blow-up}
        The map $\alpha_3$ is regular and $V_3$ is smooth.
    \end{theorem}
    \begin{proof}
        Since the union of $\widehat{{Z_1}}$ and $\widehat{{Z_2}}$ is the indeterminacy scheme of $\alpha_2$ and the center of the blow-up $\pi_3:V_3 \to V_2$, $\alpha_3$ resolves the indeterminacies of $\alpha_2$. Since $\widehat{{Z_1}}$ and $\widehat{{Z_2}}$ are blow-ups of smooth varieties at smooth centers, they are smooth. Similarly, $V_2$ is the blow up of a smooth variety along a smooth center, so $V_3$ is smooth.
    \end{proof}

\bibliographystyle{elsarticle-num} 
\bibliography{refs.bib}

\begin{thebibliography}{10}
\expandafter\ifx\csname url\endcsname\relax
  \def\url#1{\texttt{#1}}\fi
\expandafter\ifx\csname urlprefix\endcsname\relax\def\urlprefix{URL }\fi
\expandafter\ifx\csname href\endcsname\relax
  \def\href#1#2{#2} \def\path#1{#1}\fi

\bibitem{aluffiFaberPoints}
P.~Aluffi, C.~Faber, Linear orbits of d-tuples of points in $\mathbb{P}^1$, J. Reine Angew. Math. 445 (1993) 205--220.
\newblock \href {https://doi.org/10.48550/arXiv.alg-geom/9205005} {\path{doi:10.48550/arXiv.alg-geom/9205005}}.

\bibitem{aluffiFaberCurves}
P.~Aluffi, C.~Faber, Linear orbits of smooth plane curves, J. Algebr. Geom. 2 (1993) 155--184.
\newblock \href {https://doi.org/10.48550/arXiv.alg-geom/9206001} {\path{doi:10.48550/arXiv.alg-geom/9206001}}.

\bibitem{aluffiFaberAPP3}
P.~Aluffi, C.~Faber, Linear orbits of arbitrary plane curves, Mich. Math. J. 48~(1) (2000) 1--37.
\newblock \href {https://doi.org/10.1307/mmj/1030132706} {\path{doi:10.1307/mmj/1030132706}}.

\bibitem{aluffiFaberAPP}
P.~Aluffi, C.~Faber, Plane curves with small linear orbits i, Ann. Inst. Fourier (Grenoble) 50~(1) (2000) 151--196.
\newblock \href {https://doi.org/10.48550/arXiv.math/9805020} {\path{doi:10.48550/arXiv.math/9805020}}.

\bibitem{aluffiFaberAPP2}
P.~Aluffi, C.~Faber, Plane curves with small linear orbits ii, Int. J. Math. 11 (2000) 591--608.
\newblock \href {https://doi.org/10.1142/s0129167x00000301} {\path{doi:10.1142/s0129167x00000301}}.

\bibitem{tzigantchev}
D.~Tzigantchev, Predegree polynomials of plane configurations in projective space, Serdica Math. J. 34~(3) (2008) 563--596.

\bibitem{cazzadorSkauli}
E.~Cazzador, B.~Skauli, Towards the degree of the $\text{PGL}(4)$-orbit of a cubic surface, Le Matematiche 75~(2) (2020) 439--456.
\newblock \href {https://doi.org/10.4418/2020.75.2.3} {\path{doi:10.4418/2020.75.2.3}}.

\bibitem{numericalAPP}
L.~Brustenga~i Moncus\'{i}, S.~Timme, M.~Weinstein, \href{https://doi.org/10.4418/2020.75.2.2}{96120: The degree of the linear orbit of a cubic surface}, Le Matematiche 75~(2) (2020) 425--437.
\newblock \href {https://doi.org/10.4418/2020.75.2.2} {\path{doi:10.4418/2020.75.2.2}}.
\newline\urlprefix\url{https://doi.org/10.4418/2020.75.2.2}

\bibitem{equivariantAPP}
A.~Deopurkar, A.~Patel, D.~Tseng, A universal formula for counting cubic surfaces, arXiv:2109.12672 (2021).
\newblock \href {https://doi.org/10.48550/arXiv.2109.12672} {\path{doi:10.48550/arXiv.2109.12672}}.

\bibitem{fulton}
W.~Fulton, Intersection Theory, Springer-Verlag, 1984.
\newblock \href {https://doi.org/10.1007/978-1-4612-1700-8} {\path{doi:10.1007/978-1-4612-1700-8}}.

\bibitem{aluffiTensorClasses}
P.~Aluffi, Macpherson’s and fulton’s classes of hypersurfaces, Int. Math. Res. Not. (1994) 455--465\href {https://doi.org/10.1155/S1073792894000498} {\path{doi:10.1155/S1073792894000498}}.

\bibitem{aluffiCharClass}
P.~Aluffi, Computing characteristic classes of projective schemes, J. Symb. Comput. 35 (2003) 3--19.
\newblock \href {https://doi.org/10.1016/s0747-7171(02)00089-5} {\path{doi:10.1016/s0747-7171(02)00089-5}}.

\bibitem{degreeSO}
M.~Brandt, J.~Bruce, T.~Brysiewicz, R.~Krone, E.~Robeva, The degree of $\text{SO}(n,\mathbb{C})$, in: Combinatorial algebraic geometry, Vol.~80 of Fields Inst. Commun., Fields Inst. Res. Math. Sci., Toronto, ON, 2017, pp. 229--246.
\newblock \href {https://doi.org/10.48550/arXiv.1701.03200} {\path{doi:10.48550/arXiv.1701.03200}}.

\bibitem{harris}
J.~Harris, Algebraic Geometry: A First Course, Graduate Texts in Mathematics, Springer, 1992.
\newblock \href {https://doi.org/10.1007/978-1-4757-2189-8} {\path{doi:10.1007/978-1-4757-2189-8}}.

\end{thebibliography}



\end{document}